\RequirePackage{snapshot}

 \newif\ifsubsections
\subsectionsfalse

\documentclass[]{pcmi}

\usepackage[sc]{mathpazo}          
\usepackage{eulervm}               
\usepackage[scaled=0.86]{berasans} 
\usepackage[scaled=1]{inconsolata} 
\usepackage[T1]{fontenc}
\usepackage[%
	protrusion=true,
	expansion=false,
	auto=false
	]{microtype}

\usepackage{xcolor}
\usepackage{graphicx}

\ifdraft
	\definecolor{linkred}{rgb}{0.7,0.2,0.2}
	\definecolor{linkblue}{rgb}{0,0.2,0.6}
\else
	\definecolor{linkred}{rgb}{0.0,0.0,0.0}
	\definecolor{linkblue}{rgb}{0,0.0,0.0}
\fi

\PassOptionsToPackage{hyphens}{url} 
\usepackage[
    setpagesize=false,
    pagebackref,
	pdfpagelabels=false,
    pdfstartview={FitH 1000},
    bookmarksnumbered=false,
    linktoc=all,
    colorlinks=true,
    anchorcolor=black,
    menucolor=black,
    runcolor=black,
    filecolor=black,
    linkcolor=linkblue,
	citecolor=linkblue,
	urlcolor=linkred,
]{hyperref}
\usepackage[backrefs,msc-links,nobysame]{amsrefs}

\customizeamsrefs 

\theoremstyle{plain}
\newtheorem{mythm}[equation]{Theorem}
\newtheorem{myprop}[equation]{Proposition}
\newtheorem{mylem}[equation]{Lemma}

\theoremstyle{definition}
\newtheorem{mydef}[equation]{Definition}

\newcommand{\nn}{\nonumber}
\newcommand{\ms}{\medskip}
\newcommand{\msi}{\par\medskip\noindent}

\newcommand{\R}{\mathbb{R}}
\renewcommand{\H}{\mathcal H}

\newcommand{\bN}{\mathbb{N}}
\newcommand{\bS}{\mathbb {S}} %
\newcommand{\bY}{\mathbb {Y}}
\newcommand{\bH}{\mathbb {H}}
\newcommand{\bQ}{\mathbb {Q}}
\newcommand{\bT}{\mathbb {T}}
\newcommand{\bV}{\mathbb {V}}
\newcommand{\bP}{\mathbb {P}} 
\newcommand{\bD}{\mathbb {D}} %
\newcommand{\bZ}{\mathbb {Z}} %

\renewcommand{\d}{\partial}
\newcommand{\dist}{\,\mathrm{dist}}

\newcommand{\sm}{\setminus}

\newcommand{\wt}{\widetilde}

\newcommand{\ol}{\overline}
\newcommand{\ub}{\underbar}

\newcommand{\cE}{{\mathcal E}} %
\newcommand{\cN}{{\mathcal N}} %
\newcommand{\cG}{{\mathcal G}} %

\newcommand{\cF}{{\mathcal F}} %
\newcommand{\cM}{{\mathbb M}} %

\newcommand{\cQ}{{\mathfrak Q}} %
\newcommand{\cS}{{\mathcal S}} %
\newcommand{\cA}{{\mathcal A}} %

\begin{document}

%

\title{SLIDING ALMOST MINIMAL SETS AND THE PLATEAU PROBLEM}

%
%
\author{G. David}
\address{Univ Paris-Sud, Laboratoire de Math\'{e}matiques, 
UMR 8658 Orsay, F-91405 ; CNRS, Orsay, F-91405}
\email{guy.david@math.u-psud.fr}
%
%
\subjclass[2010]{Primary 49K99 ; Secondary 49Q20} 
\keywords{Almgren minimal sets, almost minimal sets, sliding boundary condition, Plateau problem}

\begin{abstract}
{
We present some old and recent regularity results concerning minimal and almost minimal
sets in domains of the Euclidean space. We concentrate on a sliding variant of
Almgren's notion of minimality, which is well suited in the context of Plateau problems
relative to soap films. We are especially interested in regularity properties near a boundary 
curve, where we would like to get a local $C^1$ description of $2$-dimensional
almost minimal sets in the spirit of J. Taylor's theorem, but we first study weaker and more
general results (local Ahlfors regularity, rectifiability, limits, monotonicity of density),
which we describe far from the boundary for simplicity.
There we insist on some simpler techniques, in particular the use of Federer-Fleming projections.
}
\end{abstract}  

%
%
\maketitle
\thispagestyle{empty}

%
%


\tableofcontents

\section{Introduction} 

The main objects of these lectures are minimal and almost minimal sets,
with the Almgren definition that seem to model best the geometry of soap films and 
bubbles. Even though this problem will not be addressed directly, we invite the reader to
think about the classical Plateau problem, where we are given a smooth closed curve 
$\Gamma$ in the Euclidean space $\R^3$, and try to find a surface (a set) spanned
by $\Gamma$ and with minimal area (Hausdorff measure of dimension $2$).
The lectures will be centered on a fairly small number of techniques, which lead
to regularity results for solutions of such a Plateau problem, and more generally
of almost minimal sets. In these notes, we shall try to insist a lot on one fundamental tool,
the so-called Federer-Fleming projection, which will be used repeatedly.
But some other constructions will be described too.
 
The archetype of a good regularity result will be the local $C^1$ description
of $2$-dimensional minimal set in $\R^3$ (far from the boundary) that was
proved by Jean Taylor in \cite{Ta1}. One of the goals of the lecture is to
describe recent attempts to give a similar description of minimal sets near  
a simple boundary of Plateau type. we will use the opportunity
given by the Park City lectures to give a simpler description of two otherwise quite
long and technical papers \cite{Sliding, C1W}, and also explain again a scheme of
proof for existence results that was introduced by V. Feuvrier \cite{Feu3}, 
and in my opinion not appreciated to its true value.

\section{Some Plateau problems} 

There are many different ways to state a Plateau problem, even in the 
most standard case when the boundary $\Gamma$ is a smooth curve in $\R^3$.
We shall mention a few in this section, but then we shall feel free to concentrate on
definitions similar to Almgren's in \cite{AlMemoir}, especially since they seem to be 
among the best models for soap films and bubbles and Joseph Plateau himself 
was interested in soap films (and also interfaces between fluids); see \cite{Plateau}.

Apparently soap films and bubbles are composed of two layers of soap molecules
with a water-repelling tail and a water-attracting head, which align themselves head 
to head with a thin layer of water in the middle; the width of the film is roughly equal 
to the length of two molecules. 
People usually come up rapidly with a simple formula for the energy of the film, 
just proportional to the total surface of the film. This sounds vague and imprecise, 
and the author finds it quite surprising that such a basic modeling actually works so well.

The general description for a solution of Plateau's problem is a set $E$
(or a similar object) spanned by $\Gamma$ and whose area is minimal,
but there are many ways to define the terms ``spanned'' and ``area'';
we shall only describe some of them.

\subsection{One dimensional sets, where Plateau is Steiner.} \label{s2.1}

The most reasonable version of Plateau's problem for $1$-d sets is the following.
Pick a finite collection of points $A_i \in \R^n$, and look for a connected set $E$ that 
contains the $A_i$ and has minimal length $\H^1(E)$ (defined below).

This is know as Steiner's problem. It is rather easy to prove that
minimizers exist, by Golab's theorem on the lower semicontinuity
of length among connected sets, and that they are composed of 
line segments whose endpoints are either points $A_i$, or some
additional points $B_j$, called \ub{Steiner points}. Near each Steiner 
point $B_j$, $E$ is composed to three line segments that end at $B_j$ 
with equal $120^\circ$ angles. This angle condition is easily proved by computing
the derivative of $\H^1(E)$ when we move $B_j$ a little and keep the other
vertices fixed. We suggest, as an exercise, to check (or at least guess) what happens 
when the $A_i$ are the four vertices of a square. Finally, $E$ has no loop
(because otherwise we may remove a segment and save some length).

Except for the invention of Steiner points and interesting questions about how fast 
one can compute the minimizers, there is not so much more to be said.
Notice however possible ruptures of symmetry (for instance, when the $A_i$ are 
the four vertices of a square), hence the lack of uniqueness,
and even more obviously the fact that some solutions are not smooth away 
from the $A_i$.

\ms
Let us give a more elaborate version of this, with nets and integer multiplicities,
which we will use as a first introduction to currents.
We are now given a finite collection of points $A_i, i\in I$, and for each one an
integer $\alpha_i \in \bZ$; we assume that 
\begin{equation} \label{2.1}
\sum_{i\in I} \alpha_i = 0,
\end{equation}
and we look for admissible nets (defined soon) with some minimality property 
that will be specified later. We see each $A_i$ as a source of electricity, 
with the intensity $\alpha_i$ (negative if if $\alpha_i < 0$), and think
of admissible nets as electrical nets that satisfy Kirchhoff's law. That is,
an admissible net is a finite collection of (oriented) intervals
$I_k = [a_k,b_k]$, $k\in K$, together with for each $k$ a multiplicity $m_k \in \bZ$,
and that satisfies the following version of Kirchhoff's law.
For each point $z\in \R^n$, denote by $K_+(z)$ the set of indices $k \in K$ such that 
$z=b_k$, by $K_-(z)$ the set of indices $k$ such that $Z=a_k$, and by
$I(z)$ the set of indices $i\in I$ such that $z=A_i$ (thus, $I(z)$ has at most one point). Then 
\begin{equation} \label{2.2}
\sum_{k\in K_+(z)} m_k - \sum_{k\in K_-(z)} m_k = \sum_{i\in I(z)}\alpha_i
\end{equation}
for all $z \in \R^n$ (but only the nodes really matter).

It is not hard to check that \eqref{2.1} is a necessary and sufficient condition on the
$\alpha_i$ for the existence of admissible nets.
We should probably require that $\alpha_i \neq 0$ for $i\in I$
(otherwise, remove $A_i$ from the discussion), that $m_k \neq 0$
(otherwise, $I_k$ is useless), and that the intervals $I_k$ have disjoint interiors 
(otherwise, we can use use a finer description where the intersection of the two
intervals is an interval of its own, with the sum of the multiplicities).

When $\cN$ is an admissible net and the intervals $I_k$ have disjoint interiors,
we shall denote by $E(\cN) = \cup_{m_k \neq 0} [a_k,b_k]$ the support of
$\cN$. Associated to $\cN$ is also the current $T = \sum_{k\in K} m_k [[I_k]]$,
where $[[I_k]]$ is a notation for the one-dimensional current of integration on 
the oriented segment $I_k$ (but we shall not define this yet).
Let us still mention that the Kirchhoff rule \eqref{2.2} is complicated way of
saying that $\d T = S$, where $S = \sum_{i\in I} \alpha_i [[A_i]]$ is the current 
of dimension $0$ associated to the data.

The simplest quantity to minimize on the class of admissible nets is probably the size
\begin{equation} \label{2.3}
S(\cN) = \sum_{k \in K ; m_k \neq 0} |b_k-a_k|,
\end{equation}
where $\cN$ denotes an admissible net, 
which is the total length of the useful part of the net. Choosing a net $\cN$ that
minimizes $S(\cN)$ corresponds roughly to minimizing $\H^1(E)$ in the Steiner 
problem above; let us just make a few observations, and leave their verification
as an exercise. 

For minimal nets, the intervals $I_k$ automatically have disjoint interiors, even
if we did not require this initially. 
If $\cN$ is a minimizer and $E = \cup_{m_k \neq 0} [a_k,b_k]$ denotes its support, 
then $E$ is as above a finite union of intervals, whose endpoints are the $A_i$ and 
Steiner points $B_j$ where exactly three intervals $I_k$ end with $120^\circ$ angles. 
Then the number of Steiner points is at most $N-2$, where $N$ is the number of 
points $A_i$, and we get a bound on the on the number of intervals $I_k$ too.
We may use this to prove that there is a minimizing net.

The set $E$ may be different from the solution of the Steiner problem above,
because it is not necessarily connected; however, we can easily choose the multiplicities
$\alpha_j$ so that any minimal net that satisfies the Kirchhoff rule has a connected support
that contains the $A_i$.
And also, if $E$ is a solution of the Steiner problem above, we can choose multiplicities
$m_i$ and $\alpha_i$ so that the associated net is supported by $E$, and even 
minimizes $S$ if all the supports of admissible nets are connected. 
This can probably be arranged, but the author was too lazy to check this; 
notice however that for a higher dimensional minimal set, the construction of a multiplicity
on this set so that the associated rectifiable current is a size minimizer is a nontrivial problem.

Instead of the size $S(\cN)$, we may also want to minimize the mass
\begin{equation} \label{2.4}
M(\cN) = \sum_{k \in K} |m_k| |b_k-a_k|,
\end{equation}
which is a more natural number when we think of $\cN$ as a current, because it is 
its norm as a linear form on the set of $1$-forms. This is the quantity that
most people like to minimize when they talk about minimal currents and surfaces.
The reader is invited to play with the size and mass minimizers that arise when
the $A_i$ are the four vertices of a square.

Let us finally mention that other quantities, such as 
\begin{equation} \label{2.4}
M_\beta(\cN) = \sum_{k \in K} |m_k|^\beta |b_k-a_k|,
\end{equation}
with $0 < \beta < 1$ are natural too, in particular in the context of
optimal networks, and produce interesting minimal nets, with angles
at the Steiner points that now depend on the $m_i$ and $\beta$.
Think about constructing an optimal net of roads that will accommodate a flux
of cars between some cities $A_i$, and where the cost of construction of
a road depends on the intensity of the traffic there, and see \cite{BCM} for information.

\subsection{Parameterizations, Rad\'{o}, and Douglas.}\label{s2.2}

Let us now think about the case when $\Gamma$ is a closed curve in $\R^3$,
and we try to minimize the area of a surface $E$ bounded by $\Gamma$, 
in the sense that $E = f(\bD)$ for some function $f$ defined on the closed unit disk $\bD$
and the restriction of $f$ to the unit circle $\bS$ is a parameterization of $\Gamma$.
The simplest way to define the area of $E$ is to (assume that $f$ is Lipschitz and)
take $\cA(f) = \int_D J_f(x) dx$, where $J_f$ is the appropriate Jacobian of $f$. 

This is what Tibor Rad\'{o} did in 1930 (see \cite{Ra1, Ra2, Ra3}), 
with conformal mappings, and Jesse Douglas in 1931, with harmonic parameterizations. 
[Of course we skip many important contributions, here and below.]
The difficulty is that the Lipschitz constant for $f$ may tend to $+\infty$ 
along a minimizing sequence, which leads to an unpleasant lack of compactness. 
They have nice solutions to this, where they first select nice parameterizations.
In particular, the paper \cite{Douglas} (which the author believes is the main reason 
for Douglas'  Fields medal) is very clever and easy to read. 

We cannot resist saying two words about it. It makes sense to decide that $f$ will be harmonic
in $B(0,1)$ and continuous on $\bD$, because such parameterizations exist. And then the area
$\cA(f)$ can be computed in terms of $f_{\vert \d S}$ alone. It turns out that the initial 
problem translates into minimizing 
\begin{equation} \label{2.6}
A(f) = \int\int \ {\sum_{j=1}^n |f_j(\theta)-f_j(\varphi)|^2
\over \sin^2\big({\theta-\varphi \over 2}\big)} \, d\theta d\varphi,
\end{equation}
where the $f_j$ are the coordinates of $f$, and this is much easier to do.

But this way to state the Plateau problem is not entirely satisfactory.
First, minimizers of $\cA(f)$ only give a good description of soap films locally when $f$
is injective. That is, if $x, y$ are interior points of $\bD$ such that $f(x) = f(y)$, then 
near $f(x)$, $E = f(\bD)$ may look like the union of two smooth surfaces that meet transversally;
soap films don't look like this, but rather like the sets of type $\bY$ that are described below.
And it is difficult to know, given $\Gamma$, when the parameterization given by Douglas will be injective.

In addition, some of the minimal sets bounded by $\Gamma$ are often best parameterized by
other sets than $\bD$, for instance with an additional handle or different topology; it is
not clear (to the author) that Douglas's argument will work in this case.

For a little more information on this and the next variants of the Plateau problem,
the reader may consult the survey \cite{SteinLecture}.

\subsection{Hausdorff measure, rectifiable sets} \label{s2.3}

From now on, all our ways to compute area will rely on the \ub{Hausdorff measure}, 
which we define now for the  convenience of the reader. 
The main properties of $\H^d$ that we like are that it is a Borel 
(but not locally $\sigma$-finite) measure, and that it coincides with 
the surface measure on smooth sets. It is defined by
\begin{equation} \label{2.7}
\H^d(E) = \lim_{\delta \to 0} \H^d_\delta(E),
\end{equation}
where 
\begin{equation} \label{2.8}
\H^d_\delta(E) = c_d \inf 
\Big\{ \sum_{j\in \bN} {\rm diam}(D_j)^d \Big\}
\end{equation}
and the infimum is taken over all coverings of $E$ by a 
countable collection $\{ D_j \}$ of sets, with 
${\rm diam}(D_j) \leq \delta$ for all $j$.

We may choose the normalizing constant $c_d$ so that $\H^d$ coincides with 
the Lebesgue measure on subsets of $\R^d$. 
See for instance \cite{Mattila} for the important verification that Borel sets
are measurable, and also information on rectifiable sets.

We shall use a lot of \ub{rectifiable sets} too. Those are the sets $E$ with 
$\sigma$-finite measure (Federer used to require finite measure, but the standard 
definitions now allow countable unions) such that $E \subset Z \cup G$, where 
$\H^d(Z) = 0$ and $G$ is a countable union of embedded $C^1$ surfaces 
of dimension $d$ (or equivalently, $G$ is a countable union of images of $\R^d$ 
by Lipschitz mappings). We will recall the properties of rectifiable sets as we use them, 
but let us already say that they have an approximate tangent $d$-plane at almost every point.

\subsection{Minimal currents}  \label{s2.4}

The most celebrated (and very successful) ways to state Plateau problems are in terms of 
currents. Existence results are made easier by stating things weakly 
(i.e., in terms of distributions), because important compactness results can be proved, 
and the setting is more in terms of differential geometry than for standard PDE's, 
relying on the integration of forms and a notion of boundary 
that comes from exterior derivatives and integration by parts. Generally, a current 
of dimension $d$ is a continuous linear form on the vector space of smooth $d$-forms 
(say, with compact support), but we will restrict here to specific classes of currents 
(rectifiable currents, integral currents) 
with a regularity which is essentially the same as the regularity of Radon measures. 
Finally regularity results for minimizers can often be proved, completing the loop and 
allowing us to return to smooth minimal sets. 
Initial and important work was done by Federer, Fleming, De Giorgi, and many others.
Out of ignorance and laziness, we just refer to \cite{Al66, Federer, FedererFleming, Fle, MoSize} 
and their references.

The simplest example of $d$-dimensional current is the current of integration on 
a smooth oriented surface $E$ of dimension $d$, which acts on a $d$-form by integrating it on 
$E$. But we are interested in the following larger classes of current, with better compactness properties.
Now let $E$ be a rectifiable set of dimension $d$, with locally finite Hausdorff measure $\H^d$ 
(defined below). Also put a measurable orientation $\tau$ (i.e., an orientation of the approximate
tangent plane to $E$ at $x$ (which exists $\H^d$-almost everywhere), chosen to be a 
measurable function of $x$),
and choose a measurable multiplicity $m(x)$ on $E$, with integer values, and integrable against
$\H^d_{\vert E}$; we define the \underbar{rectifiable current} $T$ by 
\begin{equation} \label{2.9}
\langle T,\omega \rangle = \int_E m(x) \; \omega(x)\cdot\tau(x)  \, d\H^d(x),
\end{equation}
where $\langle T,\omega \rangle$ is our notation for the effect of the current $T$ on the smooth, 
compactly supported $d$-form $\omega$, and $\omega(x)\cdot\tau(x)$ is a notation for the way 
one uses the orientation $\tau$ to integrate a form on $E$ (or a $C^1$ surface to start with).

One of the clever ideas behind the use of currents is that we can define boundaries
as in differential geometry. The \ub{boundary} $\partial T$ of the $d$-dimensional current $T$
is a current of dimension $d-1$, defined by duality by 
\begin{equation} \label{2.10}
\langle \partial T,\omega \rangle = \langle T, d\omega \rangle
\ \hbox{ for every $(d-1)$-form $\omega$,}
\end{equation}
where $d \omega$ denotes the exterior derivative of $\omega$. 
The point is that when $S'$ is the current of integration on a smooth oriented surface 
$S$ with boundary $\Gamma$, the Green formula says that $\partial S' = \Gamma'$,
the current of integration on $\Gamma$.

Notice also that  $\d\d = 0$ because $dd=0$. An \ub{integral current} is a rectifiable current $T$
(with an integrable integer multiplicity $m$) as above, such that $\d T$ is such a rectifiable current as well.

The most classical way to state the \ub{Plateau problem for currents} is to take a  
$(d-1)$-dimensional current $S$,  with $\partial S = 0$,
and minimize \ub{the mass} ${\rm Mass}(T)$, among $d$-dimensional currents $T$
that satisfy the boundary equation
\begin{equation} \label{2.11}
\partial T = S.
\end{equation}
The mass ${\rm Mass}(T)$ is the operator norm of $T$, where
we put a $L^\infty$-norm on forms. When $T$ is a rectifiable current given by \eqref{2.9},
\begin{equation} \label{2.12}
{\rm Mass}(T) = \int_E |m(x)| \, d\H^d(x).
\end{equation}
In this setting, there is a very general existence result for these 
\ub{mass minimizing} \ub{currents};
for instance, it is enough to assume that $S$ is rectifiable (as above), compactly supported,
and such that $\d S=0$ (necessary for \eqref{2.11}, because $\d\d = 0$). 
This comes from a quite strong \ub{compactness theorem}.

An in addition mass minimizing currents have good regularity properties in general.
In particular, in codimension $1$ the support of $T$ is a smooth submanifold when $n \leq 7$.  
We refer to \cite{MoBook} and its references for loads of information.

The whole theory is a great success for weak solutions and Geometric Measure Theory,
but the sad news for us is that mass minimizing currents don't describe most soap films.
For one think, there are soap films of dimension $2$ in $\R^3$ with one-dimensional
singularities, which therefore cannot come from mass minimizing currents; see the discussion
of J. Taylor's theorem below.

\subsection{Size Minimizing currents} \label{s2.5}

If we want to describe soap films, it seems that it is better to minimize the \ub{size} of $T$
among solutions of \eqref{2.11}. The size of $T$ is the Hausdorff measure of its support, i.e., 
when $T$ is given by \eqref{2.9},
\begin{equation} \label{2.13}
{\rm Size}(T) = \H^d\big(\big\{x\in E ; m(x) \neq 0\big\}\big).
\end{equation}
That is, we no longer count the multiplicity. 
The difference between mass and size minimizers is essentially the same 
as in dimension $1$ above, and for instance the minimal cones in the description
of J. Taylor's theorem can be realized as supports of size minimizers; we just need to 
find adequate multiplicities, and the minimality follows from a calibration argument.
See \cite{LaMo1, LaMo2}. 

There are some bad news about this variant of the Plateau problem. 
Some soap films are hard to describe with size minimizers, typically because 
they are not orientable; in some cases one can circumvent this problem by 
various algebraic tricks, but altogether it is a little awkward to use orientations on soap films 
which naturally are not oriented. 
See \cite{SteinLecture} for some additional information.
 
 The second piece of bad news is that there is no general existence result for size minimizers,
 even when $S$ is the current of integration on a smooth curve in $\R^3$. One cannot use the
 same compactness result as for mass minimizers, because we do not control the mass of $T$
 along a minimizing sequence. There are interesting partial results by R. Hardt and T. De Pauw
 \cite{DePauw, DePH}, with also existence results for functionals between mass and size.

Some more existence results should follow from Yangqin Fang's regularity results
from \cite{FangC1}, but only when $d=2$, $n=3$, \eqref{2.11} is replaced by a
a condition that says that $\d T$ is homologous to $S$ 
inside a smooth surface, and the support of $T$ is required to stay on one side of that surface.
 
 Another option with an effect similar to using the size is to minimize the mass, but for
 multiplicities with values in some other, often discrete group. See \cite{Annalisa}
 for an initial description.

\subsection{Reifenberg homology minimizers} \label{s2.6}

Here and below, we really consider closed sets $E$, and want to minimize $\H^d(E)$
under a topological constraint that says that $E$ is ``spanned by $\Gamma$'',
but what exactly should we mean by this?

For Reifenberg \cite{Reifenberg}, 
$E$  is a compact set that contains $\Gamma$,
and the boundary condition is stated in terms of \v{C}ech 
homology on some commutative group $G$. Specifically we require 
the inclusion $i : \Gamma \to E$ to induce 
a trivial homomorphism from $H_{d-1}(\Gamma,G)$ to $H_{d-1}(E,G)$.
Or we could instead require that it annihilates some given subgroup of $H_{d-1}(\Gamma,G)$.

When $d=2$, $n=3$, and the boundary $\Gamma$ is a curve, this is a way 
to say that $\Gamma$ (or the obvious generator of $H_{1}(\Gamma,G)$ associated 
to $\Gamma$) vanishes in $H_{1}(E,G)$, or in somewhat more vague terms, 
that $E$ ``fills the hole''.

The choice of \v{C}ech %
homology is not innocent, this homology is more complicated to define,
but it has some stability with respect to taking limits, which is good for the existence of minimizers. 
Reifenberg proved this for compact groups, such as $G = \bZ_2$ or $G = \R /\bZ$, 
under some small regularity for $\Gamma$. It is a beautiful (but tough) proof by hands, 
using minimizing sequences and haircuts.

De Pauw obtained the $2$-dimensional case when $\Gamma$ is a curve and $G = \bZ$ (with currents).
In that case the equivalence with the size minimizing problem above is not known yet (multiplicities
are hard to construct), but the infimum is the same \cite{DePauw}.

More recently, essentially optimal existence results were given by Yangqin Fang \cite{FangEx}, 
after a claim by F. Almgren with varifolds \cite{Al68}. 
A little later, following U. Menne, Y. Fang and S. Kolas\'{i}nski \cite{FaKo} 
gave a proof with flat chains inspired of Almgren's argument. 

\subsection{A sliding Plateau problem} \label{s2.7}

We finally arrive to the author's favorite definitions, based on a notion of deformations of a set.

\begin{mydef}\label{d2.14}
Let a  set $\Gamma$ be given (a priori any closed subset of $\R^n$). 
Let $E \subset \R^n$ and a closed ball $B \subset \R^n$ be given.
A \ub{sliding deformation} for $E$ in $B$, with respect to the sliding boundary $\Gamma$, 
is a one-parameter family $\{ \varphi_t \}$, $0 \leq t \leq 1$, of functions, such that
\begin{equation} \label{2.15}
(x,t) \to \varphi(x,t) = \varphi_t(x) : E \times [0,1] \to \R^n \text{ is continuous;} 
\end{equation}
\begin{equation} \label{2.16}
\varphi_0(x) = x \hbox{ for } x\in E,
\end{equation}
\begin{equation} \label{2.17}
\varphi_t(x) = x \hbox{ for $x\in E \sm B$ and } 0 \leq t \leq 1,
\end{equation}
\begin{equation} \label{2.18}
\varphi_t(E \cap B) \subset B \hbox{ for } 0 \leq t \leq 1,
\end{equation}
the important sliding condition
\begin{equation} \label{2.19}
\varphi_t(x) \in \Gamma  \ \hbox{ when } x\in E \cap \Gamma  \hbox{ and } 0 \leq t \leq 1
\end{equation}
and 
\begin{equation} \label{2.20}
\varphi_1 \hbox{ is Lipschitz.}
\end{equation}
By extension, we also call $E_1 = \varphi_1(E)$ a \ub{sliding deformation of $E$ in $B$}.
Finally, a sliding deformation of $E$ is just a set $E_1 = \varphi_1(E)$ which is a sliding deformation of 
$E$ in some $B$. 
\end{mydef}

There would be similar notions localized in an open set $U$, but we shall not need them.

This is a minor modification of a definition of Almgren \cite{AlMemoir}, to take the boundary into account.
The point is to allow $E$ to move, including along $\Gamma$, but not to be detached from $\Gamma$,
a little bit like a shower curtain.

Note that $\varphi$ is not required to be injective; in the problems below, if you can pinch $E$
and this way get a new set with less surface, this gives a good competitor.

We keep \eqref{2.20} because Almgren used it and it does not disturb, but we could also drop
it, and we should observe that no bound on the Lipschitz constant for $\varphi_1$ is ever required.

Because of the extra condition \eqref{2.19}, we really need to state things in terms of a deformation
$\{ \varphi_t \}$, rather than just the endpoint $\varphi_1$. Without \eqref{2.19} (and because
we decided to restrict to deformations in a ball; things would be different if we used deformations
in a non convex compact set), it would be easy, given the final mapping $\varphi_1$, to extend 
by convexity, set $\varphi_t(x) = t \varphi_1(x) + (1-t) x$, and observe that it provides a deformation.
But here we need to make sure that $\varphi_t(x) \in \Gamma$ for all $t$ if $x\in \Gamma$,
and when $\Gamma$ is not convex, an extension as above may be hard to find.

It would probably not be natural to demand that $\varphi_t$ be defined on the whole $B$
(as opposed to $E$ alone), because we do not really want to control the air around the soap film; 
when there is no boundary condition (or equivalently $\Gamma = \R^n$),
this makes no difference because we could extend $\varphi_t$, but in the present situation
we don't necessarily  know how to extend $\varphi_t$ so that \eqref{2.19} holds 
also for $x\in B \sm E$.

There is a Plateau problem attached to this definition: let $E_0$ be a given closed set,
and try to minimize $\H^d(E)$ among all the sliding deformations of $E_0$. 
In some cases the problem will be uninteresting, either because there is no sliding deformation 
with $\H^d(E) < +\infty$, or on the opposite the infimum is $0$. 
But in general there may be more than one interesting initial set $E_0$ for a given compact 
set $\Gamma$.

The author claims that this is probably one of the best ways to model soap films, 
and also likes the setting for the following reasons. First, the notion of sliding competitor may
fit the way soap films are created (but the author does not claim any precise knowledge about this).
It is nice that we don't need to say precisely for which topological reason $E_0$ is linked 
to $\Gamma$, or in other words to relate solutions of Plateau problems to specific topological
or algebraic reasons that may depend on the problem. 
The problem is rather insensitive to orientation. And importantly, for each
choice of $\Gamma$ there may be a few different interesting choices of $E_0$, leading to
different minimizers, just as this happens for soap films attached to a wire.

Yet there is the usual bad news: no general existence theorem is known, even when 
$d=2$, $n=3$, and $\Gamma$ is a nice curve. 

Naturally, we do not account here for unrealistic deformations that would 
extend the film too far: some real life films can be deformed into a point, 
with a long homotopy that soap is unlikely to discover because this would involve going
through a surface with a much larger area. We cannot do much about this, except for
mentioning that this may happen. Of course the dynamics of soap films and bubbles is interesting,
but this is not the subject of these lectures.

We will return extensively to the notion of sliding competitors, but in the mean time let us end 
this section with the general conclusion that many interesting Plateau problems are still there to solve.

\section{Almost minimal sets and what we want to do} \label{s3}

Here comes our last section of introduction.
Our general goal here will be to study regularity properties of sliding minimal and almost minimal
sets; in this section we give some of the relevant definitions and try to justify this goal.

Let us directly define sliding almost minimal sets in $\R^n$. Let $\Gamma \subset \R^n$
be a closed set, and consider a closed set $E \subset \R^n$ such that
\begin{equation} \label{3.1}
\H^d(E\cap B(0,R)) < +\infty \ \text{ for } R > 0.
\end{equation}
Also let $h : (0,+\infty) \to [0,+\infty]$ be a nondecreasing function such that
$\lim_{r \to 0} h(r) = 0$; we shall call this a gauge function.

\begin{mydef} \label{d3.2}
We say that $E$ is a sliding almost minimal set, with the sliding boundary $\Gamma$ 
and the gauge function $h$, when 
\begin{equation} \label{3.3}
\H^d(E \cap \ol B(x,r)) \leq \H^d(F \cap \ol B(x,r)) + h(r) r^d
\end{equation}
whenever $F$ is a sliding competitor for $E$ in any ball $B = \ol B(x,r)$ (as in Definition~\ref{d2.14}).
\end{mydef}

When we take $\Gamma = \R^n$ (i.e., forget about the sliding boundary conditions), 
we get what we'll call a \ub{plain almost minimal set}. When $h \equiv 0$ 
we get a sliding minimal set (or a plain minimal set).

It is easy to localize this definition to an open set $U$; we say that $E$ is a sliding 
almost minimal set (with the sliding boundary $\Gamma$ and gauge function $h$ 
but we shall not always repeat this) in $U$ when \eqref{3.3} holds whenever $F$ is a sliding 
competitor for $E$ in a ball $B = \ol B(x,r)$ that is contained in $U$. 
Thus Definition \ref{d3.2} corresponds to $U = \R^n$, and sliding 
almost minimal set in $U$ are automatically sliding almost minimal in any open set $U' \subset U$.
Almost all of our results will be local, i.e., concern sets that are almost minimal in (a neighborhood
of) a given ball.

Recall that one way to produce deformations of a set $E$ is to pinch them so that two pieces
of $E$ come together and we save some Hausdorff measure.
Here are some examples, starting with explicit ones, and continuing with institutional ones.

The (plain) minimal sets of dimension $1$ in a domain $U$ are composed of locally finite unions
of line segments, that can only meet by sets of three at their endpoints and with 
$120^\circ$ angles. This was checked by Morgan \cite{MoCurve} as an exercise on currents,
and later in \cite{Holder} with more elementary cut and paste argument. Then, if $U = \R^n$,
it can be checked that (modulo sets of vanishing $\H^1$ measure), the only
minimal sets are the empty set (which we'll often forget to mention), the lines, and the sets
$Y$ (three half lines with the same endpoint that make $120^\circ$ angles there)
of Subsection \ref{s2.1}.
As we have seen, the union of two transverse (or even perpendicular) lines is not 
minimal. It is not almost minimal either. 

The planes, the cones of type $\bY$ and $\bT$ described below in Subsection \ref{s7.1} 
are also plain minimal sets in $\R^n$. There are a few other explicit examples
like this (in higher dimensions), but not so many.

Then there are the minimal surfaces, like the catenoid. Those are often locally minimal only,
which means that $\H^2(E) \leq \H^2(F)$ when $F$ is a deformation of $E$ in a small enough
ball $\ol B(x,r)$. For instance, seen from far the catenoid looks a lot like two parallel planes,
that we can pinch to get a better competitor.

We should also mention that smooth surfaces (like spheres, but not only) are
rather easily seen to be locally plain almost minimal, with $h(r) \leq C r$ for $r$ small.
One of the reasons why we authorized $h(r) = +\infty$ is so that we can say without even thinking
that spheres (or objects that are even more irregular at large scales) are almost minimal.
This is also a way to make it plain that in some case we do not get any information by
comparing $E$ with a deformation in a large ball.

We shall give other examples of simple sliding minimal sets later. 
For the moment, let us only mention the case when $d=2$, $\Gamma$ is a line in $\R^3$, 
and $E$ is a half plane bounded by $\Gamma$, or the union of two half planes 
bounded by $\Gamma$ and that make an angle at of least $120^\circ$ there
(we'll call this a set of type $\bV$).

\ms
Let us turn to institutional minimal and almost minimal sets now.
If $E$ minimizes $\H^d$ among all sliding deformations of a given set $E_0$,
as in Subsection \ref{s2.7}, it is automatically sliding minimal. But there are other examples. 
First, it is often possible to prove that the limit $E_\infty$ of some minimizing sequence for the 
sliding Plateau problem above is a sliding minimal set, without being able to show that it is a 
sliding deformation of $E_0$, and studying the regularity of $E_\infty$ is useful in itself, but also
could help us prove that $E_\infty$ is a sliding deformation of $E_0$ and solve the corresponding
Plateau problem.

Also, the Reifenberg homology minimizers of Subsection \ref{s2.6} are sliding minimal sets;
see \cite{SteinLecture} for the easy verification. In \cite{FangHolder}, 
regularity properties for those in a simple context is used to prove that they are also 
solutions of the Reifenberg Plateau problem with the apparently more friendly 
singular homology. 

Similarly, the supports of size-minimizing currents of Subsection \ref{s2.5}
are sliding minimal (see \cite{SteinLecture} again)
and some regularity for those is useful in itself and could be used for existence results.

Concerning almost minimality, this is also a very useful notion, because (local) minimizers of
slight modifications of our usual functional $\H^d$ are typically almost minimal, with a gauge
function like $h(r) = Cr^{\alpha}$ (often with $\alpha=1$). The simplest example is soap bubbles,
which are also subject to different pressures from the two sides of the bubble.
The corresponding force is proportional to the surface (and the difference of pressure),
and for a deformation of $E$ in $\ol B(x,r)$, we expect a difference of potential energy of order
$r^{d+1}$, where $r^d$ comes from $\H^d(E\cap B(x,r))$ and the extra power accounts for
the displacement. For soap bubbles, we expect $E$ to have a constant mean curvature, 
proportional to the difference of pressure. Thus the pressure in smaller bubbles is 
larger, and of course with soap films the pressure is the same on both sides and the mean 
curvature vanishes. For very small soap bubbles, the pressure is very large and we expect the almost minimality constants for $E$ to deteriorate.

Of course we could include other ``small'' forces too, like the gravity, and another simple example
of almost minimal set would be a set that minimizes $\int_E f(x) d\H^d(x)$ for some 
H\"older-continuous function $f$ such that $C^{-1} \leq f \leq C$. In fact, we may also include
strongly Euclidean but yet non isotropic elliptic integrands of the form
\begin{equation} \label{3.4}
J(E) =\int_E f(x,T_E(x)) d\H^d(x),
\end{equation}
where $T_E(x)$ denotes the approximate tangent 
$d$-plane to $E$ at $x$ (we may either assume that $E$ is rectifiable, or define $J(E)$
in some other way when it is not rectifiable, and then find out after the fact that minimizers
are rectifiable), and $f$ is defined by 
\begin{equation} \label{3.5}
f(x,T) = \H^d(B(0,1) \cap A(x)(T)),
\end{equation}
where $A$ is a H\"older continuous function with 
$n \times n$-matrix values (or we should say, linear mappings of $\R^n$) such that 
$||A(x)||$ and $||A^{-1}(x)||$ are bounded. Of course the standard classes of elliptic integrands
are much larger than this (for instance, we could use $l^p$-norms instead of Euclidean ones),
but for those we should not expect the corresponding almost minimal sets to be as easy to study
as in the Euclidean case. But notice that the quasiminimal sets defined in Subsection~\ref{s9.2} 
are the same as long as $C^{-1} \leq f \leq C$.

The author likes to insist on the extra stability provided by almost minimal sets. 
While we expect roughly the same low regularity results (say, in the $C^1$ category) 
for almost minimal sets, we really like the extra flexibility. 
Other conditions, such as bounds on the first variation for varifolds, seem to be much less flexible.

\ms
So we want to study the local regularity properties of sliding almost minimal sets
(say, with gauge functions $h(r) \leq C r^{\alpha}$), possibly with the hope that 
existence results may follow. 

We will present the regularity story in two steps. First we'll describe general results
(Ahlfors regularity, rectifiability, stability under limits, blow-up limits and
minimal cones), which can also be proved near $\Gamma$ at the price of longer
and more complicated arguments, and which we shall present in the plain case
(i.e., far from $\Gamma$).

Then we will present more recent and precise results that are specific of $2$-dimensional 
sliding almost minimal sets. We think about variants of 
J. Taylor's theorem \cite{Ta1},
which will be stated in Section \ref{s7} as a best example of what we try to accomplish, 
and we'll try to describe attempts to get similar statements near points of a simple boundary 
$\Gamma$ (for instance a line). 

Finally, we shall try to explain a scheme introduced by V. Feuvrier that allows one to
use, in the circumstances where we have an appropriate regularity result for
sliding minimal sets $E$ (think about the existence of local Lipschitz retractions on $E$),
a minimizing sequence of improved competitors (so that they are quasiminimal) 
to get an existence result. This fits well with the other results presented here, 
the author's impression is that Feuvrier's argument was not well enough appreciated,
and so we want to try once more.

\section{Weak regularity properties for almost minimal sets} \label{s4}

In this section we sketch the proof of some of the weak (but general and useful)
properties of almost minimal sets, sliding or not. We rather follow the proofs of 
\cite{Sliding}, not only because they were already written in the sliding context, but
also because since time had passed, some of the initial proofs (from \cite{AlMemoir, DSMemoir},
for instance) were improved in the mean time. But we'll do the description
in the plain case, because it is less technical and otherwise almost the same.

We'll concentrate on the local Ahlfors regularity of $E$, its rectifiability, 
the rectifiability of limits, and some important stability results under limits.
The main hero for the part of the proofs that we can describe will be the Federer-Fleming
projection on dyadic cubes, also called deformation lemma in some contexts.

Our standing assumption now is that $E$ is a (coral, as in the next subsection)
almost minimal set (plain to simplify) of dimension $d$, with a small enough gauge function $h$
($h(r) \leq C r^{\alpha}$ for some $\alpha > 0$ is more than enough), in a domain
$U \subset \R^n$ which contains the balls $B(x,r)$ where we put ourselves.

Our constants $C$ will be allowed to depend on $n$, $d$, $h$, but not
on $E$, $U$, or $B(x,r)$ and the radius $r$ may have to be taken small, so that $h(r)$ 
is small enough.

\subsection{Coral (or reduced) sets.} \label{s4.1}

This will just be a precaution, so that we don't spend time discussing useless additional sets of 
vanishing $\H^d$-measure. For $E \subset U$ closed, with locally finite $\H^d$-measure in $U$,
denote by $E^\ast$ the closed support of $\H^d_{\vert E}$. That is,
\begin{equation} \label{4.1}
E^\ast = \big\{ x\in E \, ; \, \H^d(E\cap B(x,r)) > 0 
\hbox{ for all } r>0 \big\}.
\end{equation}

We say that $E$ is \ub{reduced}, or \ub{coral} when $E = E^\ast$.

If $E$ is almost minimal, then $E^\ast$ is also almost minimal,
with the same gauge $h$, because it is easy to check that $\H^d(E\sm E^\ast)=0$,
and then by direct inspection. With sliding almost minimal sets, this is also true, but it requires 
a small proof, given in \cite{Sliding}. Incidentally, when we say ``this is also true'', we mean under
quite general assumptions on the boundary sets $\Gamma$ (we can authorize more than one at the
same time, which is convenient for instance to force $E$ to lie in an initial closed domain).

So it is safe to focus on reduced sets. 
This will simplify our statements; otherwise we would typically have to say that $E$
is composed of a thin part of vanishing measure, plus a locally Ahlfors regular set, say.
But we shall keep in mind that $E^\ast \setminus E$ can play a role in some topological problems,
even if it does not show up in the nice descriptions below.
Anyway, from now on, all our sets will be coral.

\subsection{Local Ahlfors regularity.} \label{s4.2}
We just give a statement for the moment; the proof will be discussed later.
Recall the standing assumption on the coral almost minimal set $E$ in $U$.

\begin{mythm} \label{t4.2}
Local Ahlfors-regularity \cite{AlMemoir, DSMemoir}.   
\\
There exists $C \geq 1$ such that
\begin{equation} \label{4.3}
C^{-1} r^d \leq  \H^d(E \cap B(x,r)) \leq C r^d
\end{equation}
whenever
\begin{equation} \label{4.4}
x\in E \ \hbox{ and $r\in (0,1)$ are such that } B(x,2r) \subset U.
\end{equation}
\end{mythm}

\msi
Thus Ahlfors regularity is just a size condition, that says that $E$ is $d$-dimensional
in a very strong and uniform way. Here $C$ depends only on $n$ and $h$; more precisely 
we can make sure that $C$ depends only on $n$, provided that $h(r)$ is small 
enough, depending on $n$.

This sounds bland, but it is also very useful, in part because many estimates are easier to
do with Ahlfors regular sets. 
The ``hard'' part seems to be the lower bound (if $E$ is to thin, we can deform it to an 
even smaller set), but in fact both proofs use the same basic engine, the Federer-Fleming 
projections described below. 

However, in codimension $1$ the proof for the upper bound is quite simple, 
so we'll give it now. We assume that $d=n-1$, and we want to show that 
\begin{equation} \label{4.5}
\H^d(E \cap \overline B(x,r)) \leq \sigma_{n-1} r^{n-1} + h(r) r^{n-1},
\end{equation}
where $\sigma_{n-1} = \H^{n-1}(\d B(0,1))$.
As in most of our proofs, we have to find a good competitor for $E$. 
Here this will be easy. Select any point $x_0 \in B(x,r/2) \sm E$;
this is possible, since $\H^d(E \cap B(x,r/2)) < +\infty$ and $d = n-1$.
Then denote by $\pi$ the radial projection centered at $x_0$, from $B(x,r) \sm \{ x_0 \}$
to $\d B(x,r)$. Extend $\pi$ to $\R^n \sm \{ x_0 \}$ by setting $\pi(y)=y$ for
$y\in \R^n \sm B(x,r)$. Then set $\varphi_t(y) = t \pi(y) + (1-t) y$ for $y\in E$
and $0 \leq t \leq 1$; it is easy to check that $\{ \varphi_t \}$, $t\in [0,1]$, is a deformation
for $E$ in $\ol B(x,r)$. So $F = \pi(E)$ is a deformation of $E$, and the definition 
of an almost minimal set yields
$\H^d(E \cap \ol B(x,r)) \leq \H^d(F \cap \ol B(x,r)) + h(r) r^d$, as in \eqref{3.3}.
But $\H^d(F \cap \ol B(x,r)) \leq \H^d(\d B(x,r)) = \sigma_{n-1} r^{n-1}$
because $d=n-1$ and $\pi(B(x,r) \sm \{ x_0 \}) \subset \d B(x,r)$; \eqref{4.5} follows.

With a sliding boundary condition (with a nice enough boundary $\Gamma$) we would have to
be more careful (because $\pi$ does not necessarily preserve $\Gamma$). We leave the details 
as an exercise; the same sort of problems arise in higher co-dimensions, making the proof a little
less pleasant in this context.

In higher co-dimensions, and even with the apparently easier upper bound, how should we proceed?
We can still try to project $E \cap B(x,r)$ on a set with controlled $\H^d$-measure, but we'll
need to be more systematic and persistent.

\subsection{Federer-Fleming projections.} \label{s4.3}

These will be compositions of radial projections (like $\pi$ above), on faces of 
various dimensions of cubes. In some arguments, it is useful to replace cubes 
by convex polyhedra, but for the moment we stick to cubes which are much simpler to organize.
We will try to explain the construction in a friendly, but also slightly imprecise way;
the reader may consult \cite{DSMemoir}, which is not bad at all. The usual descriptions
(such as in \cite{FedererFleming} (obviously!) or \cite{Federer}) are usually a little harder 
to follow, because they are designed for currents. But the construction is the same.
We start with some notation. 

We shall use closed cubes $Q \subset \R^n$, which we could take with faces parallel to the axes 
(this costs nothing), but not necessarily dyadic to start with. 
For each such cube $Q$, $\d Q$ denotes the boundary of $Q$ in $\R^n$.
More generally we will be interested in cubes $Q \in \cQ_k$, $0 \leq k \leq n$, the set of 
$k$-dimensional (closed) cubes. Such a cube $Q$ is thus contained in some 
$k$-dimensional affine subspace $H$ of $\R^n$, and then $\d Q$ 
will denote the boundary of $Q$ in $H$. 

Notice that (for $k \geq 1$), $\d Q$ is composed of $2d$ cubes
of $\cQ_{k-1}$, which we call the faces of $Q$; we denote by $\cF_{k-1}(Q)$
the set of faces of $Q$. We also iterate, and (when $k \geq 2$) call $\cF_{k-2}(Q)$
the set of faces of cubes of $\cF_{d-1}(Q)$, and so on. 
We even call $\cF(Q)$ the union of all the $\cF_{l}(Q)$, $0 \leq l < k$; we call those 
the subfaces of $Q$.

The $l$-dimensional skeleton of $Q$ is the set $\cS_l(Q) = \cup_{S \in \cF_l(Q)}$,
a union of $l$-cubes. 

Return to the Federer-Fleming projection, and start with the main building block.
We are given a cube $Q \in \cQ_k$, contained in the $k$-dimensional affine subspace $H$,
and a point $\xi \in \frac12 Q$ (the cube of $H$ with the same center and half the sidelength);
then let $\pi_\xi$ denote the radial projection centered at $\xi$, from $Q \sm \{ \xi \}$ to $\d Q$.
Recall that $\pi_\xi(y)$ is the only point of $\d Q$ such that $y \in [\xi,y]$.
We systematically extend $\pi$ to $H \sm \{ \xi \}$ by setting $\pi(y) = y$ for
$y \in H \sm Q$; the extended mapping is still Lipschitz away from $\xi$.

In the arguments below, there is also a closed set $F \subset Q$, with $\H^d(F) < +\infty$,
and our first task is to choose $\xi \in \frac12 Q \sm F$, such that $\pi_x(F)$ is not too large. 
In fact, often we have other constraints, so we really need to show that for most choices
of $\xi \in \frac12 Q \sm F$, $\pi_x(F)$ is not too large.

There will be two estimates based on the same comment. Since $\xi \in \frac12 Q$,
all the radii $[\xi , z]$, $z\in \d Q$, are nicely transverse to $\d Q$, and there is a 
(simple geometric) constant $C$ such for each $y\in \frac12 Q \sm \{ \xi \}$, the mapping $\pi_\xi$
is $C |y-\xi|^{-1} d(Q)$-Lipschitz near $y$, where we denote by $d(Q)$ the sidelength of $Q$.
Because of this and by definition of $\H^d$ (we don't even need to disturb the area formula),
\begin{equation} \label{4.6}
\H^d(\pi_\xi(F\cap Q)) 
\leq C \int_F \Big(\frac{d(Q)}{|y-\xi|}\Big)^d d\H^d(y).
\end{equation}
The simplest estimate is most useful when we can choose $\xi$ far from $F$ as possible; 
it says that
\begin{equation} \label{4.7}
\H^d(\pi_\xi(F\cap Q)) 
\leq C  \Big(\frac{d(Q)}{\dist(\xi, F)}\Big)^d \H^d(F).
\end{equation}
We may use this sort of estimate when we have enough control on $F$, for instance if we know
that it is (a Lipschitz image of) an Ahlfors regular set of dimension $d < k$, but often we do not know
this and $F$ may be roughly dense. Fortunately, if $k$ (the dimension of $Q$) is larger than
$d$, we can use \eqref{4.6} and Fubini to show that the average value of $\H^d(\pi_\xi(F\cap Q))$
is under control, and then pick $\xi$ by Chebyshev. That is,
\begin{equation} \label{4.8}
\begin{aligned}
\int_{\xi \in \frac12 Q \sm F} \H^d(\pi_\xi(F)) d\xi
&\leq C \int_{\xi \in \frac12 Q \sm F} \int_{y \in F \sm \{ \xi \}} 
\Big(\frac{d(Q)}{|y-\xi|}\Big)^d d\H^d(y)d\xi 
\cr&\leq C \int_{y \in F} \int_{\xi \in \frac12 Q \sm \{y\}}  \Big(\frac{d(Q)}{|y-\xi|}\Big)^d
d\xi d\H^d(y)
\cr&= C d(Q)^k\int_{y \in F} \H^d(y) \leq C d(Q)^k \H^d(F)
\end{aligned}
\end{equation}
because $F$ has vanishing $k$-dimensional measure and the integral in $\xi$ converges 
when $d<k$. Thus it is easy to pick $\xi \in \frac12 Q \sm F$ such that
\begin{equation} \label{4.9}
\H^d(\pi_\xi(F)) \leq C \H^d(F).
\end{equation}

This was how we construct one basic block. Now we need to worry about gluing blocks.
We start with $k=n$. When $Q$ is a top dimensional cube, we have defined $\pi_\xi$
also on $\R^n \sm Q$, by $\pi_\xi(y) = y$. If we have a collection of cubes $Q_j$, 
disjoint except for their boundaries (and we'll only use this when they are of the same
size and belong to the same net) and for each one we pick a center $\xi_j \in \frac12 Q_j$
and define the corresponding mapping $\pi_j$, we can compose all these mappings
(in any order; they commute) and get a mapping which is equal to $\pi_j$ on $Q_j$,
and is the identity on $\R^n \sm \cup_j Q_j$.

We can also do this at the level of faces of dimension $k > d$: 
if we have a collection of faces $Q_j$ of the same dimension, disjoint except for their boundaries, 
and that belong to the same net (this will be clear when we do it; for instance, take any collection of 
$k$-faces of dyadic cubes of a given size), 
we can also compose them and get a Lipschitz mapping $\pi$ that is equal to $\pi_j$ on $Q_j$,
and the identity on the rest of the $k$-dimensional skeleton of that net.

We are now ready to iterate the basic construction and project a set of finite $\H^d$ measure
on a $d$-dimensional skeleton. For reasons that will be explained later (basically, reduce a boundary effect),
we like to do this on many small cubes at the same time. We start with more notation.

Let $Q \in \cQ_n$ be given, and let $N > 0$ be a large integer. We cut $Q$ in the obvious way 
into $N^n$ almost disjoint (i.e., with disjoint interiors) cubes $R$, $R \in \cQ(Q,N)$, of
sidelength $d(R) = N^{-1} d(Q)$, and we want to project on the faces of those cubes simultaneously.

Let us call $\cF_{k}(Q,N)$ the set of $k$-dimensional faces of those cubes (i.e., 
$\cF_{k}(Q,N) = \cup_{R \in \cQ(Q,N)} \cF_k(R)$), and also $\cS_k(Q,N)$ the corresponding
skeleton (i.e.,  $\cS_k(Q,N) = \cup_{S \in \cF_{k}(Q,N)} S = \cup_{R \in \cQ(Q,N)} \cS_k(R)$.

Our Federer-Fleming projection will be a deformation for an initial closed set $E$, 
with $\H^d(E\cap Q) < +\infty$ (think about our almost minimal set in an open set that contains $Q$).
It will be obtained by composing a collection of deformations $g_k$.

We start with $k=n$, and the set $E_n = E$. For each $R \in \cQ(Q,N)$, we apply the basic
construction above to the cube $R$ and the set $F=E_n \cap R$. We find a point
$\xi_R \in \frac12 R \sm E_n$ such that, now calling $\pi_R$ the radial projection that we called
$\pi_{\xi_R}$, 
\begin{equation} \label{4.10}
\H^d(\pi_R(E_n \cap R)) \leq C \H^d(E_n \cap R).
\end{equation}
We compose all these mappings together and find a Lipschitz mapping $g_n$, that
maps each $E_n \cap R \in \cQ(Q,N)$ to $\d R$, and is the identity on $E_n \sm Q$.
We also get that 
\begin{equation} \label{4.11}
\H^d(g_n(E_n \cap Q)) \leq C \H^d(E_n \cap Q) < +\infty,
\end{equation}
by summing \eqref{4.10} over $R$.

If $n=d+1$, we stop here.
Otherwise, we now construct $g_{n-1}$, defined on $E_{n-1} = g_n(E_n)$. 
We intend to take $g_{n-1}(y)=y$ for $y\in E_{n-1} \sm Q$, so we just need to define $g_{n-1}$
on $E_{n-1} \cap Q = g_n(E_n \cap Q)$ (because $g_n(y) = y$ outside of $Q$). 
Notice that this set is contained in $\cS_{n-1}(Q,N)$ (the union of faces of cubes $R$),
and we will define $g_{n-1}$ independently on all the faces $S \in \cF_{n-1}(Q,N)$
that compose this skeleton. For the faces $S$ that are contained in $\d Q$, we have to
take $g_{n-1}(y) = y$ on $S$, because we said we want to take $g_{n-1}(y) = y$ on $E_{n-1}\sm Q$.
For each other face $S$, we apply the basic construction to $S$, the set $F = E_{n-1} \cap S$,
and we get a point $\xi_S \in \frac12 S \sm E_{n-1}$ such that the analogue of \eqref{4.9}
holds for the radial projection $\pi_S = \pi_{\xi_S}$. Notice that $\pi_S$ coincides with the identity
on $\d S$, so that we do not have a conflict of definition on $S \cap \d Q = \d S \cap \d Q$
(since $S$ is not contained in $\d Q$). As before, we can now compose all these mappings
and get a mapping $g_{n-1}$, defined on $\cS_{n-1}(Q,N)$ minus all the centers, hence
on $E_{n-1} \cap Q$. It is not hard to check that $g_{n-1}$ is not only defined on 
$E_{n-1}$, but in fact Lipschitz there (although maybe with a large constant); there is a little
more to say, but not much, and we leave the details.

At this point we have a new set $E_{n-2} = g_{n-1}(E_{n-1})$, and the same argument as for
\eqref{4.11} also yields
\begin{equation} \label{4.12}
\H^d(E_{n-2} \cap Q)) \leq C \H^d(E_n \cap Q) < +\infty.
\end{equation}
Its image is now composed of $E \sm Q$ (where $g_n$ and $g_{n-1}$ coincide with the identity),
a piece of $\d Q$ that we were not allowed to modify, and the rest lies in the smaller
skeleton $\cS_{n-2}(Q,N)$.

If $d = n-2$, we stop. Otherwise we continue, and define $g_{n-3}$ independently
on each face $S \in \cF_{n-2}(Q,N)$. We still take $g_{n-3}(y)=y$ on $E_{n-2} \sm Q$
and on $E_{n-2} \cap \d Q$, so we take $g_{n-3}(y)=y$ on $S$ when $S \subset \d Q$.
For the other faces $S$, we apply the basic construction with $F = E_{n-2} \cap S$,
select a center $\xi_S \in \frac12 S \sm E_{n-2}$ such that the analogue of \eqref{4.9} holds,
and use the radial projection $\pi_S = \pi_{\xi_S}$. Then we define $g_{n-2}$ as before,
and continue. Eventually we get a set $E_d$, which is composed of
$E \sm Q$, a piece of $\d Q$, and a subset of the $d$-dimensional skeleton $\cS_{d}(Q,N)$.

Finally we set $f = g_n \circ g_{n-1} \ldots g_{d+1}$. It is easy to see that this is a deformation
in $Q$ (because each $g_n$ is a deformation in $Q$, for instance). This will be our standard
Federer-Fleming projection, associated to $E$, $Q$, and the large integer $N$.

It may happen that the final set $E_n = f(E)$ is so small inside of $Q$ that for each
face $S \in \cS_{d}(Q,N)$ that is not contained in $\d Q$, we can find a point
$\xi_S \in \frac12 S \sm E_{d}$. When this happens, we can continue the construction one more
step, i.e., define $g_{d}$ as above and $\wt f = g_d \circ f$, and get a new set
$E_{d-1} = \wt f(E)$ such that 
\begin{equation} \label{4.13}
E_{d-1} \cap [Q \sm \d Q] \subset \cS_{d-1}(Q,N).
\end{equation}
This is even better: we essentially managed to kill the interior of $Q$.

Let us end with two observations before we apply this to almost minimal sets.
In our construction, each mapping $g_k$ maps each $k$-face $S \in \cF_k(Q,N)$
to itself, so $g_k$, $f$, and $\wt f$ map every cube $R \in \cQ(Q,N)$ to itself.
Hence, if $V(R)$ denotes the collection all the cubes $R' \in \cQ(Q,N)$ that touch $R$
and $N(R)$ the union of these cubes,
an iteration of \eqref{4.9} yields
\begin{equation} \label{4.14}
\H^d(E_d \cap R) \leq \sum_{R' \in V(R)} \H^d(f(E \cap R')) 
\leq C \sum_{R' \in V(R)} \H^d(E \cap R') \leq C \H^d(E \cap N(R)).
\end{equation}
That is, we also control the measure of the image locally, because we know roughly where each
point of $E_d$ comes from. The same remark holds for $E_{d-1}$ when $\wt f$ is defined.

\subsection{A proof of local Ahlfors regularity} \label{s4.4}

Let us now describe a proof of Theorem~\ref{t4.2} (the local Ahlfors regularity of $E$ when 
$E$ is almost minimal).

We start with the upper bound. Let $E$ be our almost minimal set, and let
$Q$ be a cube such that $2\sqrt n Q \subset U$.
The general idea is that if $\H^d(E \cap Q)$ is too large compared to $d(Q)^d$,
the Federer-Fleming projection above gives a deformation $E_d$ of $E$ in $Q$, which in
$Q$ is essentially  contained in a $d$-dimensional skeleton, whose measure is easy to control.
A contradiction with the almost minimality of $E$ should ensue.

We do the argument with a large $N$, to be chosen later, and whose effect will be to make 
the undesirable boundary effects smaller. The difficulty comes from the contribution of
$E_d \cap \d Q$ in the estimates, which itself is a consequence of the fact that we cannot 
brutally take two definitions for $f$, a projection on $\cS_d(Q,M)$ inside $Q$ and the identity outside.

Since the almost minimality of $E$ was written in terms of balls, we use the smallest ball
$B$ that contains $Q$, whose radius is $r = \sqrt n d(Q)/2$.
Notice that $2B \subset U$ because $2\sqrt n Q \subset U$. 
We apply \eqref{3.3} with this $B$, remove the contribution
of $B \sm Q$ which is the same for $E$ and $F = E_d$, and get that
\begin{equation} \label{4.15}
\H^d(E \cap Q) \leq \H^d(E_d \cap Q) + h(r) r^d.
\end{equation}
Now write $E_d \cap Q = F_1 \cup F_2$, where $F_1 = E_d \cap \d Q$,
and hence $F_2 \subset \cS_d(Q,N)$ by construction. The contribution
of this part is easily estimated, since
\begin{equation} \label{4.16}
\H^d(F_2) \leq \sum_{R \in \cQ(Q,N)} \H^d(\cS_d(R))
\leq C N^n (N^{-1} d(Q))^d = C N^{n-d} d(Q)^d.
\end{equation}
This looks large because of $N$, but recall that we can pick the Ahlfors regularity constant
after we choose $N$. For $F_2$, denote by $\cQ^1(Q,N)$ the set of cubes $R \in \cQ(Q,N)$ 
that touch $\d Q$, and set $A^2(Q,N)= \cup_{R \in \cQ^1(Q,N)} N(R)$; this is a thin
annulus in $Q$ near $\d Q$. We observe that $F_1 \subset f(A^2(Q,N))$, apply 
\eqref{4.14}  to each cube $R \in \cQ^1(Q,N)$, notice that the sets $N(R)$ have bounded overlap, 
and get that
\begin{equation} \label{4.17}
\begin{aligned}
\H^d(F_1) &\leq C \sum_{R \in \cQ^1(Q,N)} H^d(F_1 \cap R)
\leq C \sum_{R \in \cQ^1(Q,N)} H^d(E\cap N(R)) 
\cr&\leq C H^d(E\cap A^2(Q,N)).
\end{aligned}
\end{equation}
Thus by \eqref{4.15}-\eqref{4.17},
\begin{equation} \label{4.18}
\begin{aligned}
\H^d(E \cap Q) &\leq \H^d(E_d \cap Q) + h(r) r^d
\cr&\leq C H^d(E\cap A^2(Q,N)) + C N^{n-d} d(Q)^d +  h(r) r^d.
\end{aligned}
\end{equation}
Let us assume that $d(Q)$, and then $r$ are so small that $h(r) \leq 1$, say, and rewrite
\eqref{4.18} as
\begin{equation} \label{4.19}
\H^d(E \cap Q) \leq C H^d(E\cap A^2(Q,N)) + C(N) r^d.
\end{equation}
If $\H^d(E \cap Q) \leq 2C(N) r^d$, we are happy because we get an upper bound
on $\H^d(E \cap Q)$. Otherwise, we get the information that
\begin{equation} \label{4.20}
H^d(E\cap A^2(Q,N)) \geq (2C)^{-1} \H^d(E \cap Q),
\end{equation}
which is strange because $A^2(Q,N)$ is as thin as we want, so it should not
bring such a large contribution to $\H^d(E \cap Q)$.

The standard way to continue the argument (see for instance \cite{DSMemoir})
would be to apply the argument again to a cube $Q_1$ slightly larger than $Q$, 
so that $Q = Q_1 \sm A^2(Q_1,N)$), find that $H^d(E\cap A^2(Q_1,N))$ is even larger, 
iterate with larger cubes $Q_k$ and values $N_k$ of $N$ that get larger each time, 
and eventually find that if $N$ was chosen large enough, all the $Q_j$ are contained 
in $2Q$ and then $\H^d(E \cap 2Q) = +\infty$.

We do not give the details here because they are a little painful and, once we get to \eqref{4.20}, 
the significant part of the argument is actually done. For the construction of an appropriate sequence 
of cubes, we refer to Lemma 4.3 in \cite{Sliding}, starting near (4.23).

\ms
This completes the proof of the upper bound in \eqref{4.3}; for the lower bound,
we actually proceed in a similar way. Again we start with a cube $Q$ such that $\sqrt n Q \subset U$, 
assume that $d(Q)^{-d} \H^d(E \cap Q)$ is very small, and try to reach a contradiction.

Let us again apply the Federer-Fleming argument to $Q$ with the large integer $N$
(to be chosen later). This gives a set $E_d$, which in the interior of $Q$ is contained
in the skeleton $\cS_d(Q,N)$. For each $d$-cube $S \in \cF_d(Q,N)$ of that skeleton,
if $S$ is not contained in $\d Q$, we can apply \eqref{4.14} to $S$ and find that
\begin{equation} \label{4.21}
\H^d(E_d \cap S) \leq C \H^d(E\cap N(S)) \leq C \H^d(E\cap Q) < \frac12 \H^d(\frac12 S)
\end{equation}
if $d(Q)^{-d} \H^d(E \cap Q)$ is small enough (depending on $N$). That is, we are in 
the case when we can apply one more projection, define $\wt f$, and use the set $E_{d-1}$ 
which is also a deformation of $E$ in $Q$ (and hence in the smallest ball $B$ that contains $Q$).

We compute as before, but now $F_2$ is contained in a $(d-1)$-dimensional skeleton, so
we get that 
\begin{equation} \label{4.22}
\H^d(E \cap Q) \leq C H^d(E\cap A^2(Q,N)) + h(r) r^d
\end{equation}
instead of \eqref{4.18}, where $r = \sqrt n d(Q)/2$ as before. Again $C$ does not depend on $N$.

This is suspicious, because on average $H^d(E\cap A^2(Q,N))$ should be much smaller
than $\H^d(E \cap Q)$, and then \eqref{4.22} would say that $\H^d(E \cap Q) \leq C h(r) r^d$,
which is even smaller than expected.
The standard way to proceed would be (as in \cite{DSMemoir}) to iterate the argument, find 
a sequence of cubes $Q_j$ that are concentric with $Q$ and whose density tends to $0$,
and then observe that this is not possible if $Q$ is centered on a point of positive upper density.

Here is a hint on how we can proceed otherwise (following the argument below 
Lemma 4.39 in \cite{Sliding}). 
If the lower Ahlfors regularity fails, we can find $Q$ as above such that 
\begin{equation} \label{4.23}
\H^d(E\cap Q) \leq c d(Q)^d,
\end{equation}
with $c$ and $d(Q)$ as small as we want, but also
\begin{equation} \label{4.24}
\H^d(E\cap \frac12 Q) \geq  c 2^{-d} d(Q)^d.
\end{equation}
For this we assume that the center of $Q$ is a point of upper density $1$ for $E$,
and replace $Q$ by $\frac12 Q$,  $\frac14 Q$, and so on, until \eqref{4.23} fails for
the first time.

Now $\H^d(E\cap Q) \leq 2^d \H^d(E\cap \frac12 Q)$.
We use this and Chebyshev (in fact, the pigeon hole principle) to replace $Q$ 
by a concentric cube $Q_1$, with $\frac23 Q \subset Q_1 \leq Q$, and for which 
\begin{equation} \label{4.25}
\H^d(E \cap \wt A^2(Q_1,N)) \leq C N^{-1} \H^d(E \cap Q_1) \leq C N^{-1} \H^d(E \cap Q).
\end{equation}
We still have \eqref{4.22} for $Q_1$, so 
\begin{equation} \label{4.26}
\begin{aligned}
\H^d(E\cap \frac12 Q) &\leq \H^d(E \cap Q_1) \leq C H^d(E\cap A^2(Q_1,N)) + h(r) r^d
\cr&\leq C N^{-1} \H^d(E \cap Q) + h(r) r^d,
\end{aligned}
\end{equation}
with a constant $C$ that does not depend on $N$,
and this contradicts \eqref{4.23} or \eqref{4.24} if $h(r)$ and $N^{-1}$ are small enough.
\qed

\subsection{rectifiability, uniform rectifiability, and projections} \label{s4.5}

Recall the definition of rectifiability in Subsection \ref{s2.3}. 
It was proved by Almgren \cite{AlMemoir} that plain almost minimal sets are rectifiable. 
In fact, plain almost minimal sets are even uniformly rectifiable (UR) \cite{DSMemoir}. 
We do not want to say too much about this, but let us at least give a statement 
with the relevant definition.
We still assume that $E$ is a coral almost minimal set in $U$ with gauge function $h$.

\begin{mythm} \label{t4.27}
Uniform rectifiability with BPLG \cite{DSMemoir}.  
\\
There exists $\theta > 0$ and $M \geq 0$ such that
for each choice of $x\in E$ and $r\in (0,1)$ such that $B(x,2r) \subset U$,
there is a $d$-dimensional Lipschitz graph $G$, with Lipschitz constant at most $M$,
such that 
\begin{equation} \label{4.28}
\H^d(E \cap B(x,r) \cap G) \geq \theta r^d.
\end{equation}
\end{mythm}

The Lipschitz part means that there is a $d$-plane $P \subset \R^n$ and an
$M$-Lipschitz function $\psi : P \to P^\perp$ such that
$G = \big\{ y+\psi(y)\, ; \, y\in P \big\}$. 

Provided that $h(r) \leq 1$, say, the constants $M$ and $\theta$ depend only on $n$
and $d$.

The property stated in the theorem  ($E$ contains big pieces of Lipschitz graphs locally) 
is in fact the combination of two properties. First, $E$ is locally uniformly rectifiable, 
which has many equivalent definitions (see \cite{Asterisque , UR}). 
A simple one is that $E$ locally contains big pieces of Lipschitz images of balls in $R^d$. 
That is, there exist $\theta > 0$ and $M \geq 0$
such that, for $x\in E$ and $r\in (0,1)$ such that $B(x,2r) \subset U$, we can find an 
$M$-Lipschitz function $g: \R^d \cap B(0,r) \to \R^n$ such that 
\begin{equation} \label{4.29}
\H^d(E \cap B(x,r) \cap g(\R^d \cap B(0,r))) \geq \theta r^d.
\end{equation}
This one is obviously weaker, but in fact not that much. If it is satisfied and in addition 
$E$ has big projections locally, then it contains big pieces of Lipschitz graphs locally.
Big projections mean that we can find $\theta >0$ such that, for $x\in E$ and $r\in (0,1)$ 
such that $B(x,2r) \subset U$, we can find a $d$-plane $P \subset \R^n$ such that
\begin{equation} \label{4.30}
\H^d(\pi(E \cap B(x,r)) \geq \theta r^d,
\end{equation}
where $\pi$ denotes the orthogonal projection onto $P$. The converse (big projections imply
BPLG) is clear. See \cite{Asterisque , UR}.

We will see a proof of the rectifiability of $E$ that also works for sliding almost minimal sets
(as soon as $\Gamma$ is reasonable, and with a little more work). Surprisingly, even though
the almost minimality is itself a quantitative notion, the uniform rectifiability of $E$ is really
complicated to get, and the author does not know how to extend it to sliding almost minimal sets,
except in simple cases where the uniform rectifiability near the boundary is not really significant
because it follows too easily from the same thing far from the boundary.
It is also worth noticing that in terms of proving other results, uniform rectifiability is not
as indispensable as the author once believed.

Let us now prove that
\begin{equation} \label{4.31}
\text{every almost minimal set is rectifiable,}
\end{equation}
in a way that can be extended to sliding almost minimal sets. The basic tool will be,
once again, Federer-Fleming projections. We start with the observation, that the author owes to
V. Feuvrier \cite{FeuThesis} (but may have been known from Almgren), 
that if $F$ is totally unrectifiable, $Q$ is cube of dimension $k > d$, and 
$\pi_\xi$ denotes the radial projection on $\d Q$ with the center $\xi$ 
(as in the early part of Subsection \ref{s4.3}), then
\begin{equation} \label{4.32}
\H^d(\pi_\xi(F \cap {\rm int}(Q)) = 0
\ \text{ for $\H^{k}$-almost every } \xi \in \frac12 Q. 
\end{equation}
Recall that we say that the set $F$ is totally unrectifiable when $\H^d(F \cap G) = 0$
for every rectifiable set $G$ (or equivalently, for every $C^1$ embedded submanifold of dimension $d$).

This looks like the Besicovitch-Federer 
Projection Theorem, except that the projections are not parallel to $(n-d)$-planes,
and in fact the proof of \eqref{4.32} uses that theorem (and Fubini).

So let us prove \eqref{4.31}. Assume instead that the almost minimal set $E$ is not rectifiable.
Write $E = E_{rect} \cup E_{irr}$, where $E_{rect}$ is rectifiable and $E_{irr}$ is totally unrectifiable
(see \cite{Mattila} for this and the density properties below).
Since $\H^d(E_{irr}) >0$, we can find $x\in E$ such that the upper density of $E_{rect}$
at $x$ vanishes. That is,
\begin{equation} \label{4.33}
\lim_{r \to 0} r^{-d} \H^d(E_{rect} \cap B(x,r)) = 0.
\end{equation}

Now consider a small cube $Q$ centered at $x$, and perform a Federer-Fleming projection as before,
except that when we choose the centers $\xi_S$ for the various faces $S$, we use the standard
Chebyshev argument to make sure that the following things happen. 
Let $k$ denote the dimension of $S$; thus $k \leq n$ and we are in the middle of the 
construction of $g_k$. 
The set $E_k$ has a part in $\d S$, which will not change because $\pi_\xi(y) = y$
for $y\in \d S$. The totally unrectifiable part of $E_k \cap {\rm int}(S)$ is sent to a negligible set;
this can easily be done because of \eqref{4.32}. Finally, the $\H^d$-measure of the 
rectifiable part of $E_k \cap {\rm int}(S)$ is multiplied by at most $C$; this last can be arranged by
Chebyshev and the proof of \eqref{4.9}. Since we know now that $E$ is locally Ahlfors regular,
we could even choose $\xi_S$ far from $E_k$, so that $\pi_\xi$ is $C$-Lipschitz and the argument
looks a little bit simpler, but let us not bother yet.

When we proceed like this, $\H^d$-almost all of the unrectifiable part of 
$E \cap {\rm int}(Q) \sm \cS_d(Q,N)$ disappears, because it is contained in the union 
of the interiors of the subfaces $S$ of dimensions $k > d$ that are not contained in $\d Q$.
The unrectifiable part of $E \cap \cS_d(Q,N)$ is negligible, because 
the faces of dimension $d$ are $d$-rectifiable, so they do not really meet $E_{irr}$.

We are left with the rectifiable part of $E \cap Q$. 
If $Q$ is small enough, $\H^d(E_{rec} \cap Q) \leq \varepsilon d(Q)^d$, with $\varepsilon$ 
as small as we want, because its center $x$ was chosen so that \eqref{4.33} holds.
Then applying $f$ multiplies this measure by at most $C$, by choice of the $\xi_S$.
Altogether,
\begin{equation} \label{4.34}
\H^d(f(E\cap {\rm int}(Q))) \leq C \varepsilon d(Q)^d.
\end{equation}
This is small enough for $f(E)$ not to contain the full $\frac12 S$ for any face $S \in \cF_d(Q,N)$.
Then we can proceed as we did below \eqref{4.21}, compose $f$ with a last mapping $g_{d}$
onto a set $E_{d-1}$ which in the interior of $Q$ is contained in a $(d-1)$-dimensional skeleton.
Then \eqref{4.22} holds, and we may conclude as above, or more simply observe that since we now
know that $E$ is locally Ahlfors regular, we could easily have used Chebyshev to choose
a cube $Q'$, such that $\frac12 Q \subset Q' \subset Q$, and for which \eqref{4.22} fails 
because $H^d(E\cap A^2(Q',N)) \leq C N^{-1} \H^d(E \cap Q) \leq C N^{-1} d(Q)^d$
and $\H^d(E \cap Q') \geq \H^d(E \cap \frac12 Q) \geq C^{-1} d(Q)^d$. 
The desired contradiction comes from applying the argument above to $Q'$.
\qed

\ms
We end this section with a remark on how to find big projections. We claim that 
when $E$ is flat, it has no big hole. Here is the corresponding statement.

\begin{mythm} \label{t4.35}
Keep $E$ as above. For each $\tau \leq 0$ we can find $\varepsilon > 0$ and 
$r_0 > 0$ so that the following holds.
Let $x\in E$ and $r \in (0,r_0]$ be such that $B(x,2r) \subset U$. 
Let $P$ be a $d$-plane through $x$ and suppose that
\begin{equation} \label{4.36}
\dist(y,P) \leq \varepsilon r \ \text{ for } y \in E \cap B(x,r).
\end{equation}
Let $\pi_P$ denote the orthogonal projection on $P$. Then
\begin{equation} \label{4.37}
\pi_P(E \cap B(x,r)) \supset P \cap B(x,(1-\tau)r).
\end{equation}
\end{mythm}

This stays true (with appropriate modifications) in the sliding case.
It applies in many small balls because $E$ is rectifiable and locally Ahlfors Regular, 
so it has approximate tangent planes almost everywhere (see \cite{Mattila}), 
and these approximate tangent planes are actual tangent planes 
(see Exercise 41.21 in \cite{MSBook}). 

Once we know that $E$ is locally uniformly rectifiable, this also implies that it 
has big projections (and then that it contains BPLG), because the local uniform rectifiability
gives enough balls where $E$ is well approximated by $d$-planes (look for the 
WGL in \cite{Asterisque}).

\ms
Let us first give a proof when $d = n-1$.
Let $E$, $B(x,r)$, and $P$ satisfy the assumptions.
Set $B' = B(x,(1-\tau)r)$.
Suppose we can find $\xi \in P \cap B' \sm \pi_P(E \cap B(x,r))$. 
We want to deform $E$ in $\ol B(x,r)$ and save some area.

Draw the vertical line $\pi_P^{-1}(\xi)$ that does not meet $E \cap B(x,r)$.

Move points of $E \cap B(x,r)$ parallel to $P$, away from the line $\pi_P^{-1}(\xi)$, and
send them radially to the thin vertical wall 
$W = \big\{ z\in \pi_P^{-1}(P \cap \d B' \, ; \, \dist(z,P) \leq \varepsilon r \big\}$.

The measure in the tube of the deformation $F$ is at most 
$\H^d(W) \leq C \varepsilon r^d$.

All the measure $\H^d(E \cap B(x,r/2)) \geq C^{-1} r^d$ disappeared. 

This contradicts the almost minimality if $h(r)$ is small enough. 

The proof in higher co-dimension is not much harder. 
We first project $E \cap B \cap \pi_P^{-1}(B')$ on $P \cap B'$, 
and do a nice interpolation of the mapping between that set and $\d B(x,r)$,
where we want our deformation to be the identity. We get a mapping which is nearly 
$1$-Lipschitz, because $E$ stays so close to $P$ in $B(x,r)$. Then only, once the points
are sent to $P$, we push them radially in $P \cap B'$, starting from the center
$\xi$ which is still not in the image. This way $W = P \cap \d B'$ costs nothing, and the
necessary gluing inside of $B(x,r) \sm \pi_P^{-1}(B')$ costs as little as we want.
See Lemma~10.10 of \cite{DSMemoir} for the proof (of almost the same statement) and 
Lemma 7.38 or 9.14 in \cite{Sliding} for the more complicated sliding version.
\qed

\section{Limits of almost minimal sets} \label{s5}

There are many things that we would expect to be true, but have to be proved.
The main one is the following.

\begin{mythm} \cite{Limits, Sliding} \label{t5.1}
Let $U \subset \R^n$ be given, and suppose that each of the sets
$E_k$ is a coral (sliding), almost minimal set in $E$, always with the same
reasonably nice sliding boundary $\Gamma$ and the same gauge function $h$.
Suppose in addition that $\{ E_k \}$ converges, locally in $U$ to a closed
set $E_\infty$. Then $E_\infty$ is a coral (sliding), almost minimal set in $E$, 
with the same sliding boundary the same gauge function $h$.
\end{mythm}

Reasonably nice allows $\Gamma$ to be a $C^1$ surface of any dimension,
but more complicated choices are allowed.

We will define convergence very soon, but there will be no surprise.

Theorem \ref{t5.1} extends to quasiminimal sets (defined in
Subsection \ref{s9.2} below). 

There is also a variant where the sliding boundary for $E_k$ is a set
$\Gamma_k$ that converges nicely to $\Gamma$, but let us skip it for the 
moment. 

When $E_k$ is almost minimal with a gauge function $h_k$, and
$\lim_{k \to +\infty} h_k(r) \equiv 0$, then $E_\infty$ is in fact minimal;
this follows rather easily, because the statement shows that it is minimal
with any of the functions $h_N = \sup_{k \geq N} h_k$.

In fact, there is even a statement that says that as long as the $E_k$
are almost minimal with a fixed gauge function, or even quasiminimal
with a fixed constant, and in addition it is a minimizing sequence, then 
$E_\infty$ is minimal. We will use this for \eqref{9.9} in Section \ref{s9.2}. 

\ms
Theorem \ref{t5.1} will be the main topic of this section.
In the present state of affairs, it is still too complicated to be entirely explained
here, but at least we will be able to say something about the lower semicontinuity
estimate that is at the center of the argument. 
Maybe soon Camille Labourie will come up with a simpler argument
for the full limiting theorem.
First we define the convergence, using the following normalized
local Hausdorff distance: for $x, r$ such that $B(x,r) \subset U$, set
\begin{equation}\label{5.2}
\begin{aligned}
d_{x,r}(E,F) &= r^{-1} \sup\big\{ {\rm dist}(y,E) \, ; \, 
y\in F \cap B(x,r) \big\}
\cr& \hskip 1.2cm 
+ r^{-1} \sup\big\{ {\rm dist}(y,F) \, ; \, 
y\in E \cap B(x,r) \big\}.
\end{aligned}
\end{equation}
When $F \cap B(x,r)$ is empty, we let the first supremum be $0$, and similarly for 
the second supremum when $E \cap B(x,r) = \emptyset$.
Then let the $E_k$ be closed in $U$  (we don't need the other case) and 
$E_\infty$ be closed in $U$ (this gives the simplest definitions); 
we say that $\{ E_k \}$ tends to $E_\infty$ (locally in $U$) when
\begin{equation}\label{5.3}
\lim_{k \to +\infty} d_{x,r}(E_k, E) = 0 \ 
\hbox{ for every ball $B(x,r) \subset \subset U$.}
\end{equation}
This definition is nice, because a standard argument with diagonal subsequences shows that
given any sequence $\{ E_k \}$ of closed sets in $U$, we can always extract a subsequence
that converges (locally in $U$) to some closed set $E_\infty$.

Our main example will be the \ub{blow-up limits} of a closed set $E$.
Given $x_0 \in E$, a blow-up limit of $E$ at $x_0$ is any limit of a convergent
sequence $\{ E_k \}$ (as above), where $E_k = r_k^{-1} [E - x_0]$ for some
sequence $\{ r_k \}$ that tends to $0$. By what we just said, there is always at least
one blow-up limit of $E$ at $x_0$ (start from $r_k = 2^{-k}$ and extract a converging subsequence),
and in general there may be lots of blow-up limits of $E$ at $x_0$ (think about a spiral).

When we work with sliding boundaries, we typically take sequences for which in addition the dilations 
$r_k^{-1} [\Gamma - x_0]$ converge to a limit $\Gamma_\infty$.

So we decided to study the limits of a sequence $\{ E_k \}$ of almost minimal sets
in $U$, all with the same $U$, the same boundary set $\Gamma$ (to simplify),
and the same gauge function $h$. An important ingredient is the following.

\begin{mylem}  \label{t5.4}
Let $U$, and $\{ E_k \}$ be as above, and suppose that $\{ E_k \}$
converges to $E_\infty$. Then $E_\infty$ is rectifiable.
\end{mylem}

The author believes that Almgren \cite{AlMemoir} probably had this (in the plain case) with
essentially the proof below, but he did not check recently. In \cite{DSMemoir}, the authors proved first
that the sets $E_k$ are uniformly rectifiable (with uniform bounds), and deduced the lemma 
from this; this looks reasonable, because uniform rectifiability goes to the limit well, while
simple rectifiability does not. That is, it is very easy to find a sequence of rectifiable sets $E_k$
(for instance composed of $4^k$ little squares) that converge to a totally unrectifiable Cantor set
$E_\infty$. But of course the sets $E_k$ are not uniformly almost minimal!

The author only re-discovered Lemma \ref{t5.4} (with probably stupid surprise) after
spending some time, not being able to prove the uniform rectifiability of the
sliding almost minimal sets and being upset about it. In addition, the proof is quite simple
and it seems that rectifiability of the limit is nearly as useful as uniform rectifiability.

So let us give (the idea of) the proof in \cite{Sliding}, except that in order
to simplify the argument a little, we forget about the sliding boundary and assume that
the $E_k$ are plain almost minimal sets.

We suppose that $E_\infty$ is not rectifiable, and proceed as in the proof of rectifiability.
Take a point $x_0 \in E_\infty$, where the upper density of the rectifiable part $E_{\infty}^{rec}$
vanishes. Then let $Q_0$ be any small cube centered at $x_0$. 

First observe that the $E_k$ satisfy the Ahlfors regularity properties in $3Q_0$, 
uniformly in $k$. A simple covering argument (with balls of the same size) shows that
then $E_\infty$ also is Ahlfors regular in $2Q_0$. 

Next we can replace $Q_0$ by a concentric cube $Q'$, such that 
$\frac12 Q \subset Q' \leq Q$, and for which the measure of $E_\infty$ near $\d Q$ is 
fairly small. More precisely, cut $Q$ into $N^n$ subcubes $R \in \cQ(N,Q)$ as we did before, 
and set $A_+(N,Q) = (1+N^{-1})Q' \sm (1-10N^{-1})Q'$, designed to be a little larger than 
the thin annulus $A^2(N,Q)$ that was used before; by the pigeon hole principle, we 
we can easily find $Q$ such that 
\begin{equation} \label{5.5}
\H^d(E_\infty \cap A_+(N,Q)) \leq C N^{-1} \H^d(E_\infty \cap 2Q_0) 
\leq C N^{-1} d(Q_0)^{d}.
\end{equation}
Next we can use again the local Ahlfors regularity of the $E_k$ and $E_\infty$,
and simple coverings by balls, to prove that for $k$ large,
\begin{equation} \label{5.6}
\H^d(E_k \cap A^2(N,Q)) \leq C \H^d(E_\infty \cap A_+(N,Q)) 
\leq C N^{-1} d(Q_0)^{d}.
\end{equation}

Let us now find a Federer-Fleming projection $f$ that essentially kills 
$E_{\infty}^{irr} \cap Q$, and multiplies the (very small) measure of $E_{\infty}^{rec}$ 
by at most $C$. We proceed a little differently as for \eqref{4.31}, because it will
be better to have some uniformity. So, when we choose the centers $\xi_S$ 
of the various faces $S$ in the Federer-Fleming construction, we 
use the local Ahlfors regularity of $E_\infty$ near $2Q$ to select $\xi_S$
at distance at least $C^{-1}N^{-1}d(Q)$ from the previous image of $E_\infty$.
We skip the details again, but the reader may find this argument in 
Lemma 3.31 of \cite{DSMemoir}, and later in \cite{Sliding}.

This way, \eqref{4.7} says that all the mappings that compose $f$ are $C$-Lipschitz. 
Of course $f$ is $C$-Lipschitz too. Notice that it is also naturally defined and
still $C$-Lipschitz (as a composition of radial projections) in a small neighborhood
of $E_\infty$, which in particular contains $E_k \cap 2Q$ for $k$ large.
Then $f$ never multiplies the measure of pieces of $E_k$ inside cubes by more than $C$, 
and we get that for $k$ large,
\begin{equation} \label{5.7}
\H^d(f(E_k) \cap A(N,Q)) \leq C \H^d(E_k) \cap A^2(N,Q))
\leq C N^{-1} d(Q_0)^{d},
\end{equation}
where $A(N,Q)$ is the thinner annulus composed of cubes $R \in \cQ(N,Q)$ that touch
$\d Q$, and by \eqref{5.6}. This takes care of $f(E_k) \cap A(N,Q)$. 

As in the proof of \eqref{4.31}, we still have some latitude to choose the centers $\xi_S$,
in particular so that in $Q \sm A(N,Q)$, the image of $E_\infty^{irr}$ is negligible.
Since $d(Q)^{-d} \H^d(E_\infty^{rec} \cap Q)$ is as small as we want, this
implies that $f(E_\infty) \cap Q \sm A(N,Q)$ never fills a half face. This allows us to 
add one more step to the construction, and get a new mapping, which we shall
also call $f$, so that now $f(E_\infty) \cap (Q \sm A(N,Q))$ is contained in a 
$(d-1)$-dimensional skeleton, and (again by construction of $f$), 
this is also true for $f(E_k) \cap (Q \sm A(N,Q))$ for $k$ large.
Together with \eqref{5.7}, this yields
\begin{equation} \label{5.8}
\H^d(f(E_k) \cap Q) \leq C N^{-1} d(Q_0)^{d}.
\end{equation}
Recall that $\H^d(E_k \cap D) \geq C^{-1} d(Q_0)^{d}$ by local Ahlfors regularity.
When $Q$ is sufficiently small, so that $h(\sqrt n d(Q))$ is very small, all this contradicts 
the almost minimality of $E_k$; the rectifiability of $E_\infty$ follows.
\qed

\ms
We come to the main reason why Theorem \ref{t5.1} works, which is the 
lower semicontinuity of $\H^d$.

\begin{mylem}  \label{t5.9}
Let $U$, and $\{ E_k \}$ be as above, and suppose that $\{ E_k \}$
converges to $E_\infty$. Then 
\begin{equation} \label{5.10}
\H^d(E_\infty \cap V) \leq  \liminf_{k \to +\infty} \H^d(E_k\cap V)
\end{equation}
for every open set $V \subset U$.
\end{mylem}

This even works for quasiminimal sets (see Subsection \ref{s9.2}), 
and also when we replace $\H^d$ with a large class of \ub{elliptic integrands} \cite{FangEx}. 
Here we'll rapidly discuss the case of plain almost minimal sets, but sliding minimal 
(or even sliding quasiminimal) sets work as well.

For the lower semicontinuity property, it makes sense to take $V$ open. For a closed square $V$,
for instance, it could be that the $E_k$ are lines segments outside $V$ that tend to a
side of $V$.  It is worth noting that the result fails miserably without the almost minimality 
assumption: a sequence of dotted lines may converge to a line. But dotted lines are 
not almost minimal!

In earlier versions of \cite{Limits} and a first part of \cite{Sliding}, the author
insisted on using the fact that our almost minimal sets satisfy the 
``uniform concentration property'' of Dal Maso, Morel, Solimini \cite{DMS}.
This property was introduced, in the context of minimizers of the Mumford-Shah functional,
precisely to prove the lower semicontinuity of $\H^{n-1}$ along minimizing sequences,
and then possibly get existence results, and it was very tempting to use it, especially
because it is an easy consequence of uniform rectifiability. It turns out \cite{Sliding} to be also a 
consequence of the rectifiability of limits (Lemma \ref{t5.4}), which is lucky because
we still cannot prove yet that sliding almost minimal sets are always uniformly rectifiable.
But in fact Yangqin Fang \cite{FangEx} discovered that there is a simpler direct proof
of \eqref{5.10}, that also uses the rectifiability of the limit, and works in the context
of elliptic integrands. 

Let us say a few words about the proof. 
Let $\{ E_k \}$, $E_\infty$, and $V$ be as in the statement.
Let $\varepsilon > 0$ and $\tau > 0$ be small, to be chosen later.
Observe that since $E_\infty$ is rectifiable, for $\H^d$-almost every point 
of $E_\infty \cap V$, we have that for $r > 0$ small enough,
\begin{equation} \label{5.11}
\text{$E_\infty$ is $\varepsilon r$-close to some $d$-plane $P = P(x,r)$ in $B(x,2r)$,}
\end{equation}
which just means that $\dist(y,P) \leq \varepsilon r$ for $y\in E \cap B(x,2r)$, 
and also
\begin{equation} \label{5.12}
\H^d(E_\infty \cap B(x,r)) \leq (1+\tau) \omega_d r_i^d,
\end{equation}
where $\omega_d$ is the $\H^d$-measure of the unit ball in $\R^d$.

The balls that satisfy \eqref{5.11} and \eqref{5.12} are what is often called a Vitali
covering of $E \cap V$, and by Vitali's covering argument (see for instance \cite{Mattila})
we can cover $\H^d$-almost all of $E_\infty \cap V$ by disjoint balls $B(x_i,r_i) \subset V$,
with small radii $r_i$, and that satisfy \eqref{5.11} and \eqref{5.12}. Then we can find
a finite subcollection $\{ B(x_i,r_i) \}$, $i\in I$, that catches most of the mass, so that
\begin{equation} \label{5.13}
\H^d(E_\infty \cap V) \leq \varepsilon + \sum_{i\in I} (1+\tau) \omega_d r_i^d.
\end{equation}
Now we use Theorem \ref{t4.35} to estimate each $r_i^d$.
Notice that there is a finite number of indices $i$ to try, and for each one, 
if $k$ is large enough the assumptions of Theorem \ref{t4.35} are satisfied
by $E_k$, with a ball $B_{i,k}$ centered on $E_k$ (this is needed for the statement),
contained in $B(x_i, r_i)$ (this will be used soon), and yet of radius $(1-10\varepsilon) r_i$ 
(this is easy to obtain). The size assumption on the radius $r$ of the theorem is satisfied 
if we made sure to take the radii small enough in our Vitali collection earlier. 
And maybe \eqref{4.36} is only satisfied with the constant $2\varepsilon$. 
Anyway, if we choose $\varepsilon$ small enough, depending on $\tau$, we get that
\begin{equation} \label{5.14}
\begin{aligned}
\H^d(E_k \cap B(x_i, r_i)) &\geq \H^d(E_k \cap B_{i,k})
\geq \H^d(\pi_P(E_k \cap B_{i,k})) 
\cr&\geq [(1-\tau)(1-10\varepsilon)]^d \omega_d r_i^{d},
\end{aligned}
\end{equation}
where $P$ and the last estimate come from \eqref{4.37}.
We sum all this, use the fact that the $B(x_i,r_i)$ are disjoint and contained in $V$, 
and get that for $k$ large,
\begin{equation} \label{5.15}
\begin{aligned}
\H^d(E_\infty \cap V) &\leq \varepsilon + \sum_{i\in I} (1+\tau) \omega_d r_i^d
\cr&
\leq \varepsilon + (1+\tau) [(1-\tau)(1-10\varepsilon)]^{-d}
 \sum_{i\in I} \H^d(E_k \cap B(x_i, r_i))
 \cr&
 \leq \varepsilon + (1+\tau) [(1-\tau)(1-10\varepsilon)]^{-d} \H^d(E_k \cap V).
\end{aligned}
\end{equation}
Since $\varepsilon$ and $\tau$ can be chosen as small as we want, \eqref{5.10}
and Lemma \ref{t5.9} follow.
\qed

\ms
At this point, it really feels like Theorem \ref{t5.1} should be easy to prove.
As these notes are being written, this is not the case. Both in \cite{Limits}
(for the plain case) and \cite{Sliding} (for the sliding case), an additional long and painful
construction of competitors is used with lots of special cases and coverings.
But the author hopes that Camille Labourie will soon come up 
with a much more pleasant (and even more general) proof soon.

The conclusion of this section is that, thanks in particular to the (surprising)
rectifiability of limits, we have a very nice tool, Theorem \ref{t5.1}, that will allow
us to use compactness in many circumstances and make it simpler to think about
regularity. The story below about blow-up limits and minimal cones depends on this!

The author imagines that a large part of people's preference for weak objects such as currents,
flat chains, varifolds, was largely due to the apparent absence of a good limiting theorem for sets.
Hopefully Theorem \ref{t5.1} will also become easy in the near future.

\section{Monotonicity of density, near monotonicity, and blow-up limits} \label{s6}

Here we start with plain almost minimal sets; the situation for sliding almost minimal sets
is more complicated and will be discussed later.

It is very useful to know, in many circumstances involving minimality (think about
minimal sets, surfaces, and currents, but free boundary problems are also concerned by this issue),
that some scale-invariant quantity is nondecreasing, or just nearly monotone.
Here the quantity of interest will be the density
\begin{equation} \label{6.1}
\theta_x(r) = r^{-d} \H^d(E \cap B(x, r)),
\end{equation}
where we often take the origin $x$ in the (almost) minimal set $E$.
We shall learn with time that low density often rhymes  with simplicity, so saying that
$\theta_x$ is nondecreasing can also be a way to say that the situation in smaller balls
tends to be simpler.

The monotonicity of $\theta_x$ for minimal things is a rather ubiquitous fact. 
Here we just give a small number of statements and hints of proofs, then discuss
variants, easy applications, and the sliding case.

\subsection{Near monotonicity in the plain case.} \label{s6.1}

We start with the simplest statement:
\begin{equation} \label{6.2}
\begin{aligned}
&\text{If $E$ is a plain minimal set in $U$, then for any $x\in U$}
\cr&\ \ \text{ $r \to \theta_x(r)$ is nondecreasing on $(0,\dist(x,\d U))$.}
\end{aligned}
\end{equation}

But we shall often use the more general near monotonicity of $\theta_x$
for almost minimal sets $x$. Even though some statements (depending on which
definition of almost minimality one takes) also work for $x\in U \sm E$, 
we shall restrict our attention to $x\in E$.

\begin{mythm} \label{t6.3}
There is $\alpha > 0$, that depends on $n$, such that 
if $E$ is a plain almost minimal set in $U$, with a gauge function $h$ that satisfies 
a Dini condition, then for $x\in E$
\begin{equation} \label{6.4}
 r \to \theta_x(r) \, 
{\rm exp \,} \Big\{ \alpha \int_0^r h(2t) {dt \over t} \Big\}
\end{equation}
is a nondecreasing function on $(0,\frac12\dist(x,\d U))$.
\end{mythm}

By Dini condition, we just mean that the integral $\int_0^r h(2t) {dt \over t}$ 
converges near $0$. When $h \equiv 0$ and $x\in E$, we recover \eqref{6.2}.
Near monotonicity is almost as useful as monotonicity, because the Dini condition says that
the integral in \eqref{6.4} has a limit when $r$ tends to $0$, so that for 
instance, under the assumptions of the theorem, we can define the density at 
$x \in E$ by 
\begin{equation} \label{6.5}
\theta(x) = \theta_x(0) = \lim_{r \to 0} \theta_x(r).
\end{equation}
Notice that $C^{-1} \leq \theta(x) \leq C$ by Theorem \ref{t4.2}; in fact it is not
hard to see that $\theta(x) \geq \omega_d$ (the density of a plane).

Let us discuss the proof of Theorem \ref{t6.3}. 
We start with the main case when $E$ is minimal, as in \eqref{6.2}, 
and the main idea is to compare $E$ with the cone over $E \cap \d B(x,r)$.
Let us assume for simplicity that $x=0$.

\begin{figure}[!h]  
\centering
\includegraphics[width=10.cm]{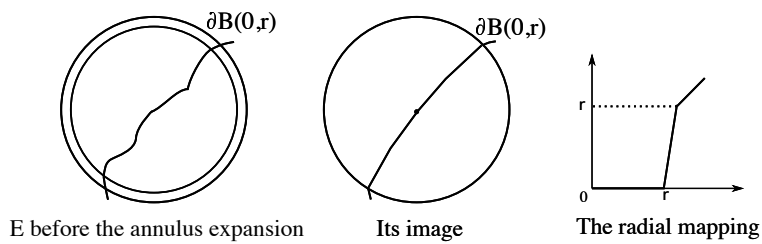}
\caption{A deformation of $E$ that is close to the cone. 
\label{fcone}}
\end{figure}

Since $r \to \H^d(E\cap B(0,r))$ is nondecreasing, 
it is the integral of its derivative (seen as a Stiljes measure),
which is no less than its almost-everywhere derivative. 
Thus, after a computation that we skip here, 
it is enough to check that for almost every $r \in (0,\dist(0,\d U))$,
\begin{equation} \label{6.6}
r^{-d}\,{\partial \over \partial r} \big(\H^d(E\cap B(0,r)\big)
\geq d \, r^{-d-1} \H^d(E \cap B(0,r)).
\end{equation}
But, by the co-area theorem, or rather more directly by approximating the rectifiable
set $E$ by $C^1$ surfaces and computing the derivative for each one,
\begin{equation} \label{6.7}
{\partial \over \partial r} \big(\H^d(E\cap B(0,r)\big)
\geq \H^{d-1}(E \cap \partial B(0,r))
\end{equation}
for almost every $r$. So it is enough to show that for a.e. $r$, 
\begin{equation} \label{6.8}
\H^d(E \cap B(0,r)) \leq {r \over d} \  \H^{d-1}(E \cap \partial B(0,r)).
\end{equation}
This is beginning to look good, because if $X$ denotes the 
the cone over $E \cap \partial B(0,r))$, another application of the co-area formula
(or more simply a parameterization of the cone and the area formula) yields
\begin{equation} \label{6.9}
{r \over d} \,  \H^{d-1}(E \cap \partial B(0,r)) = \H^d(X \cap B(0,r)).
\end{equation}
So the proof would be finished if we knew that $X$ coincides in $\ol B(0,r)$
with a deformation of $E$ in $B(0,r)$. This is not exactly the case, but yet 
we can approach $X$ (in $\ol B(0,r)$) by Lipschitz deformations of $E$ in $B(0,r)$,
where the Lipschitz mapping is radial, expands a lot an annulus near $\partial B(0,r)$,
and contract the rest of $B(0,r)$ brutally to the origin. See Figure\ref{fcone}

Thus, with a simple limiting argument it is easy to get \eqref{6.2}.
In the case of almost minimal sets, the proof is a little longer but the idea remains the same.
As before, for almost each $r$, we can approximate the cone $X$, construct
deformations of $E$ in $\ol B(0,r)$, apply the definition of almost minimal sets, and get
an estimate that relates $\H^d(E\cap B(0,r))$ to $\H^{d-1}(E\cap \d B(0,r))$.
We can see this as a differential inequality that relates $\H^d(E\cap B(0,r))$
and its distribution (or Stiljes) derivative, which we can integrate to recover \eqref{6.4}.
That is, if we differentiate the right-hand side of \eqref{6.4}, the derivative that comes from 
the Dini integral compensates the error term that comes from the almost minimality
of $E$, if $\alpha$ is large enough. The computations are slightly more unpleasant 
than in the minimal case, but there is no surprise at the end.
\qed

\subsection{The almost-constant density principle} \label{s6.2}

Here we describe a trick which is often useful, and relies on (unfortunately more than)
the proof of Theorem \ref{t6.3}. In rough terms, if $E$ is almost minimal with a small enough
gauge $h$, and in addition $\theta_x$ is nearly constant (instead of nearly monotone)
on $(0,2r)$, then $E$ looks a lot like a minimal cone in $B(x,r)$.
This turns out to be quite helpful.

Of course a preliminary to this is that if $E$ is a minimal set in $B(0,r)$,
and $\theta_0$ is constant on $(0,r)$, then $E$ is a minimal cone 
(we should say, coincides with a minimal cone in $B(0,r)$).

This is true, but the author only knows a surprisingly unpleasant proof \cite{Holder}.
When we look at the proof above, we rapidly find that for almost every $x\in E$,
the tangent plane to $E$ at $x$ contains the origin. But to go from this to the 
conclusion, the author did not find a way that does not use the construction
of a complicated deformation of $E$.

Let us now state the almost-constant density principle, and then say how it follows 
from the preliminary fact.

\begin{myprop} \label{t6.10}
For each small $\delta > 0$, we can find $\varepsilon >0$
such that, if $x\in E$ (an almost minimal set in $U$, with gauge function $h$), 
$B(x,2r) \subset U$, $h(2r) \leq \varepsilon$, and 
\begin{equation} \label{6.11}
\theta_x(2r) \leq \inf_{0 < t < 10^{-3}r}\theta_x(t) + \varepsilon,
\end{equation}
then there is a minimal cone $Z$ centered at $x$ such that
\begin{equation} \label{6.12}
d_{x,r}(E,Z) \leq \delta
\end{equation}
and even (a form of approximation in measure)
\begin{equation} \label{6.13}
\big| \H^d(E\cap B(y,t))-\H^d(Z\cap B(y,t)) \big|
\leq \delta r^d \ \ \text{ for } B(y,t) \subset B(x,r).
\end{equation}
\end{myprop}

The infimum in \eqref{6.11} is a little strange; we put it for technical reasons 
(to avoid requiring a Dini condition on $h$), but you could think about $\theta_x(0)$
instead. Notice that in \eqref{6.13} the normalization is by $r^d$, not $t^d$, so \eqref{6.13} 
is not interesting for $t$ small. Although we did not say, $\varepsilon$ depends on
$n$, but nothing else. See Proposition 7.24 in \cite{Holder} in the plain case, and
Proposition 30.19 in \cite{Sliding} otherwise.

The proof use compactness and limits in a standard way.
Suppose the proposition fails for some $\delta$, and let $E_k$ provide a counterexample
for $\varepsilon_k = 2^{-k}$, in some ball $B_k$. By translation and dilation invariance 
of the problem, we may assume that $B_k = B(0,1)$. And the gauge function $h_k$
associated to $E_k$ is such that $h_k(2) \leq \varepsilon_k = 2^{-k}$

Then use the compactness result evoked below \eqref{5.3} to extract a subsequence 
(that will still be denoted by $\{ E_k \}$) such that $\{ E_k \}$ converges to a limit 
$E_\infty$ locally in $B(0,2)$.
By Theorem \ref{t5.1}, the discussion below its statement about sequences $\{ h_k \}$ 
that tend to $0$, and the fact that we never need a ball larger than $B(0,2)$ in the arguments, 
$E_\infty$ is minimal in $B(0,2)$. 

The next step is to start from \eqref{6.11} for $E_k$ and take a limit.
We know from Lemma \ref{t5.9} that 
\begin{equation} \label{6.14}
\H^d(E_\infty \cap B(0,2)) \leq \liminf_{k \to +\infty} \H^d(E_k \cap B(0,2))
\end{equation}
Then by \eqref{6.11} $\H^d(E_k \cap B(0,2)) 
\leq 2^{d-k} + 2^d \inf_{0 < t < 10^{-3}}[t^{-d}\H^d(E_k \cap B(0,t))]$; we fix any
small $t$ and get that 
\begin{equation} \label{6.15}
\H^d(E_k \cap B(0,2)) \leq 2^{d-k} + 2^d t^{-d}\H^d(E_k \cap \ol B(0,t)).
\end{equation}
Now we need a result of upper semicontinuity, that says that 
\begin{equation} \label{6.16}
\limsup_{k \to +\infty} \H^d(E_k \cap \ol B(0,t)) \leq \H^d(E_\infty \cap \ol B(0,t)).
\end{equation}
Let us just say a few words about the proof of this, and refer to
Lemma 3.12 of \cite{Holder} in the plain case or Theorem 21.1 in \cite{Sliding} otherwise.
The main ingredient is again the rectifiability of the limit $E_\infty$, which together with
a covering argument allows us to reduce to a good upper bound for 
$\H^d(E_k \cap B(x,r))$ when $E_\infty$ (and then $E_k$ for $k$ large) is very close 
to a $d$-plane $P$. As usual, we build a deformation, which in this case is quite close
(in $B(x,r)$) to the orthogonal projection on $P$. This allows us to compare 
$\H^d(E_k \cap B(x,r))$ with $\H^d(P \cap B(x,r))$, with small errors that come from
gluing near $\d B(x,r)$.

So we get \eqref{6.16}, we compare with the previous estimates, find that 
$t^{-d}\H^d(E_\infty \cap \ol B(0,t)) \geq 2^{-d}\H^d(E_\infty \cap B(0,2))$
for $t$ small, use the monotonicity of the density $\theta(t) = t^{-d}\H^d(E_\infty \cap B(0,t))$
(because $E_\infty$ is minimal) to show that $\theta(t)$ is constant on $(0,1)$, and then
our preliminary fact to show that $E_\infty$ is a cone. 

We then get the desired contradiction almost in the expected way. Since $E_\infty$
is the limit of the $E_k$, \eqref{6.12} holds for $k$ large, and so does \eqref{6.13},
but this time with the help of Lemma \ref{t5.9} (for the lower semicontinuity)
and the proof of \eqref{6.16} (for the upper semicontinuity). 
At some point we need to be a little careful about the measure of spheres, but nothing bad
ever happens.
\qed

\subsection{Blow-up limits} \label{s6.3}

Recall from a few lines below \eqref{5.3} that a blow-up limit of $E$ at $x$ is any limit
of a sequence $\{ E_k \}$, where $E_k = r_k^{-1} [E - x]$ for some 
sequence $\{ r_k \}$ that tends to $0$. We expect that in some cases, 
$E$ may have more than one blow-up limit at $x$, but the author only knows 
easy counterexamples (spirals) when the gauge function $h$ decays very slowly.

Let $x\in E$ be a given plain almost minimal set. 
We shall assume that $h$ satisfies a Dini condition, so that 
$\theta(x) = \theta_x(0) = \lim_{r \to 0} \theta_x(r)$ is defined, as in \eqref{6.5}.
We claim that for each blow-up limit $X$ of $E$ at $x$,
\begin{equation} \label{6.17}
X \text{ is a minimal cone in $\R^n$,}
\end{equation}
with the (constant) density $\theta(x)$, i.e., 
\begin{equation} \label{6.18}
\H^d(X\cap B(0,R)) = \theta(x) R^d \ \text{ for } R > 0.
\end{equation}
The fact that $X$ is minimal comes from Theorem \ref{t5.1},
because it is easy to see that $E_k$ is almost minimal in the large open set 
$r_k^{-1}(U-x)$, with the small gauge function $h(r_k \cdot)$.
Then 
\begin{eqnarray} \label{6.19}
R^{-d} \H^d(X\cap B(0,R)) &\leq& R^{-d}\liminf_{k \to +\infty} \H^d(E_k \cap B(0,R))
\nn \\ &=& \liminf_{k \to +\infty} (r_k R)^{-d}\H^d(E \cap B(x,r_k R))
= \lim_{r \to 0} \theta_x(r) = \theta(x)
\end{eqnarray}
by algebra and \eqref{6.5}. For the other inequality, we use the upper semicontinuity
estimate \eqref{6.16} and get that 
\begin{eqnarray} \label{6.20}
R^{-d} \H^d(X\cap B(0,R)) &=& R^{-d} \H^d(X\cap \ol B(0,R))
\geq R^{-d}\limsup_{k \to +\infty} \H^d(E_k \cap \ol B(0,R))
\nn \\ &=& \limsup_{k \to +\infty} (r_k R)^{-d}\H^d(E \cap \ol B(x,r_k R))
=\lim_{r \to 0} \theta_x(r) = \theta(x);
\end{eqnarray}
\eqref{6.18} follows.

Usually minimal cones of dimension $d$ are somewhat easier to study than
general minimal sets of dimension $d$, because when $X$ is a minimal cone,
$K = X \cap \d B(0,1)$ is a $(d-1)$-dimensional object (thus a priori much simpler),
and it has some almost minimality properties in $\d B(0,1)$ that it inherits
from the almost minimality of $X$. So it makes sense to decide to study first
the minimal cones $X$, and then try to deduce local regularity properties
for $E$ near $x$ from the knowledge of its blow-up limits at $X$.
We play this game a lot in the next sections.

\subsection{Monotonicity of density for sliding almost minimal sets?} \label{s6.4}

What about the sliding almost minimal sets then?
Let us first assume, to simplify the discussion, that we study the density 
$\theta_0(r)$ at the origin.

In the proof of monotonicity, we compare $E$ with the cone $X$ 
over of $E \cap \d B(0,r)$, or approximations of $X$. This involves mappings that
move the points along radial direction, and in general, these mappings will not satisfy the
sliding condition \eqref{2.19} and our proof will not extend.

Neither does the (near) monotonicity of density. The simplest example is a half plane $H$
bounded by a line $\Gamma$. It is not so hard to show (and extremely easy to imagine) that
$H$ is a sliding minimal set of dimension $2$, with sliding boundary $\Gamma$. 
Suppose that $0 \in H$ and $\dist(0,\Gamma) > 0$. Then $\theta(0,r)=\pi$ for 
$r < \dist(0,\Gamma)$, then it decreases, and eventually tends to $\pi/2$ at $+\infty$; 
this is not what we wanted.

Yet there is a simple case when the arguments given above all work well, which is when 
\begin{equation} \label{6.21}
\text{$\Gamma$ is a cone and we take $x \in \Gamma $.} 
\end{equation}
In this case of balls centered on the cone $\Gamma$, the results mentioned 
above still work for sliding almost minimal sets, with no important modification.

We can also extend them to the similar case when $0 \in \Gamma$ and 
$\Gamma$ is almost like a cone, for instance a smooth manifold, and we get a small additional term
that can be incorporated in the almost monotonicity formula \eqref{6.4}.

Unfortunately it is also very good to get near monotonicity results for balls
centered away from $\Gamma$; we will very rapidly describe in Subsection \ref{s8.2}
some partial results with a different functional \cite{Mono}, and how we can try to use them.

\section{Minimal cones and the Jean Taylor theorem} \label{s7}

We start with the program described in Subection \ref{s6.3}; since the blow-up limits
of almost minimal sets are minimal cones, why not try to list all the minimal cones, and
then deduce regularity results from this? This program works extremely well for plain almost minimal
sets of dimension $2$ in $\R^3$; this is the theorem of Jean Taylor \cite{Ta2} that we 
describe in this section.

\subsection{Minimal cones of dimension $1$ or $2$.} \label{s7.1}

We want to make lists of minimal cones.
We start with dimension $1$. Recall that we promised not to mention the empty set. 
Then the \ub{minimal cones of dimension $1$} in $\R^n$ are just the lines and the $Y$-sets, 
i.e. the unions of $3$ half lines that end at the same point, 
where they make equal angles of $2\pi/3$. 
This is quite easy, because we just have to consider the set $K = E \cap \d B(0,1)$, 
which is finite, observe that if the cone over $K$ is to be minimal, the angle of any two 
of the half lines through $K$ that emanate from $0$ has to be at least $2\pi/3$, 
because otherwise we can pinch them near the origin and get a better competitor, 
and finally that if $3$ points of $\d B(0,1)$ ever make an angle at least $2\pi/3$ 
with each other, they have to lie in a same plane and make angles of $2\pi/3$.

It is just a bit more complicated to check that the list of minimal sets of dimension $1$
in (the whole) $\R^n$ is the same. We leave it as an exercise with a hint: you may 
prove this by hand, but it may be more fun to use as much as you can of the theory
of the previous sections, including limits and the monotonicity of density. 
Otherwise see \cite{MoCurve}. 

We now switch to the \ub{minimal cones of dimension $2$ in $\R^3$}.
The following list was given a long time ago by Ernest Lamarle (a student of J. Plateau);
the proof was completed later by A. Heppes and then 
J. Taylor. See \cite{Lamarle, Heppes, Ta2}. 

\begin{figure}[!h]  
\centering
\includegraphics[width=6.cm]{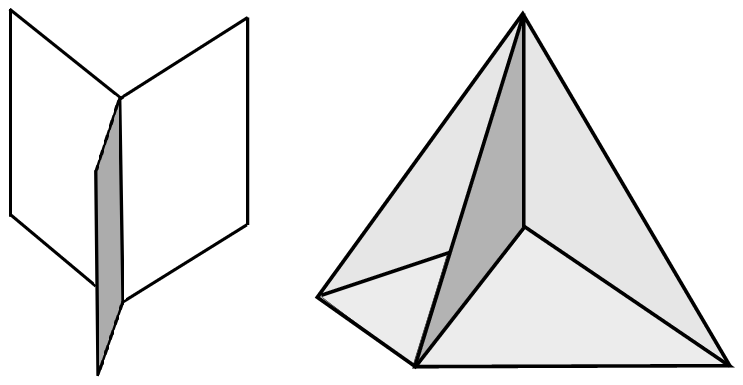}
\caption{Left: a cone of type $\bY$. Right: a cone of type $\bT$. \label{fYT}}
\end{figure}

These are the \ub{planes} (we may also say sets of type $\bP$),
the \ub{cones of type $\bY$}, which are the unions of three half planes bounded by a
same line $L$, and that make angles of $2\pi/3$ along $L$,
and the \ub{cones of type $\Bbb T$}
(see Figure \ref{fYT}). The latter are obtained by drawing 
a regular tetrahedron centered at the origin, and then taking the (positive) cone 
over the union of the $6$ edges of the tetrahedron; this gives a cone with $6$ faces that touch 
each other at the origin and otherwise meet with each other with angles of $2\pi/3$.

Let us just say a few words about the proof.
The fact that the cones of type $\bY$ (easy by slicing too) $\bT$ are minimal is usually 
proved with a calibration argument; see \cite{LaMo1}.
Then we need to check that no other cone $X$ is minimal.
It is not too hard to give a combinatorial description of $K = X \cap \d B(0,1)$:
it is a finite union of arcs of geodesics (great circles or arcs of great circles), 
that may only meet at their extremities and with angles of $2\pi/3$. 
The unpleasant part of the proof then consists in making a full list of cones for which $K$
has the description above, and then finding a better competitor for each of those, except
the three above.
See Ken Brakke's home page \url{http://www.susqu.edu/brakke/} 
for pictures of the minimal cones, the unlucky candidates 
(such as the cone over the edges of a cube, or a prism), 
and a better competitor for each one.

Unfortunately, we do not have a full list of minimal cones of dimension $2$ in
$\R^n$, $n \geq 4$. We still have the combinatorial description of 
$K = X \cap \d B(0,1)$ as a finite union of great circles or arcs of great circles that
meet with $2\pi/3$ angles, and we know a few more examples, but even when $n=4$
we expect that our list is incomplete. Let us merely mention that 
Xiangyu Liang \cite{Liang2P, Liang 3P} 
was able to prove that the almost orthogonal union of two (or more if $\R^n$ contains them)
planes is minimal (the orthogonal case was known, with a calibration argument), 
and that the product of two sets $Y$ that lie in orthogonal planes is minimal \cite{Liang2P}.

Obviously, minimal cones of dimension $d=3$ or higher will be even more challenging.
See \cite{LuuHolder, Luu3D} for very first attempts in this direction, and the more recent impressive 
\cite{CES}.

\subsection{Jean Taylor's theorem at last.} \label{s7.2}

Let $E$ be a plain coral almost minimal set of dimension $2$ in $U \subset \R^3$,
with a gauge function $h$ such that $h(r) \leq C r^{\alpha}$ for some $\alpha > 0$.
Then let $x\in E$ be given. We know that $E$ has at least a blow-up limit at $x$,
and that for every such blow-up limit $X$, $X$ is a minimal cone such that  
$\H^d(X\cap B(0,1)) = \theta(x)$ as in \eqref{6.18}.

This leaves us with three options: $\theta(x) = \pi$ and all the blow-up limits of $E$ at $x$
are planes; $\theta(x) = 3\pi/2$ and all the blow-up limits of $E$ at $x$ are
sets of type $\bY$, and $\theta(x) = \theta_T$, the density of cones of type $\bT$,
and all the blow-up limits of $E$ at $x$ are sets of type $\bT$.
In either case, it turns out that $E$ is $C^{1+a}$-equivalent to any of its  
blow-up limits of $E$ at $X$ in a neighborhood of $x$.
We give a slightly more precise statement, and then comment.

\begin{mythm} \label{t7.1}
Let $E$ be an almost minimal set of dimension $2$ in $U \subset \R^3$, 
with a gauge function $h(r) \leq C r^{\alpha}$ for some $\alpha > 0$. 
Let $x\in E$ be given, and let $Z$ be a blow-up limit of $E$ at $x$.
Then there is a radius $r > 0$ and $C^{1+a}$-diffeomorphism 
$\Phi : B(0,2r) \to \Phi(B(x,2r)) \subset \R^3$ such that
$\Phi(0) = x$ and $\Phi(Z \cap B(0,2r)) = E \cap \Phi(B(0,2r)) \supset E \cap B(x,r)$.
\end{mythm}

Here $a>0$ is a constant that depends only on $n$ and $\alpha$. 
Notice that this description implies the uniqueness of the blow-up
limit of $E$ at $x$ (in this case, we'll say $E$ has a tangent cone).

We can even say more. For each small $\tau > 0$, we can make sure that
$(1-\tau)|y-z| \leq |\Phi(y)-\Phi(z)| \leq (1-\tau)|y-z|$ for $y, z \in B(0,2r)$,
and then naturally $\Phi(B(0,2r)) \supset B(0,r)$.

It is important that not only $E$ has a $C^{1+a}$ parameterization, but also 
the parameterization extends to the ambient space. But in the present case, all this
simply boils down to the fact that $E$ is composed of the right number of faces,
that meet with the correct $2\pi/3$ angles. The story about $2\pi/3$ angles
was not announced in the statement, but it easily follows because the blow-up
limits of $E$ at such points must be sets of type $\bY$ (except naturally for the central
point when $Z$ is of type $\bT$).
 
More regularity (even, all the way to analyticity) 
could be obtained if we were concentrating
on minimal sets. Yet we don't know the best value of $a$ (but the author would expect $\alpha/2$).
In the other direction, we expect that some (Dini-type) constraint on $h$ is needed for the theorem
to hold, we know how to make it work with gauge functions significantly larger than $r^\alpha$,
but did not really look for the optimal results.

This is a really beautiful result. Recall that the singularities above can be seen
in real soap films. The theorem says that (as long as the modeling is correct) you
won't see any other ones. This is what we would love to imitate in more complicated situations.

Yet let us mention one drawback: at this time, we do not always know how to estimate $r$. 
This is related to the following unpleasant fact.
Suppose $E$ is almost minimal in $B(0,1)$, with a small enough gauge function, and is close
enough in $B(0,1)$ to a plane or a set of type $\bY$. Then our proof shows that the 
conclusion of Theorem \ref{t7.1} holds, and that we can even take $r = 1/2$. 
But when $E$ looks a lot like a cone of type $\bT$, we cannot say that $E$ contains
a point of type $\bT$ (i.e., where the blow-up limits of $E$ are cones of type $\bT$).

For {$2$-dimensional sets in $\R^n$}, $n \geq 4$, there is a similar result \cite{Holder, C1}
but weaker in the sense that unfortunately we do not have a precise complete list of minimal cones,
and we can prove the full $C^{1+a}$ result only for some blow-up limits; for the other ones, we only 
get a H\"older equivalence and maybe not the uniqueness of the blow-up limit.

The author tried in \cite{Montreal} to give a shorter description of the proof, and will use this excuse to
say very little here. Surprisingly, for the H\"older equivalence, the main ingredients are merely
the monotonicity of density, the almost-constant density principle of Subsection \ref{s6.2} 
and (to find points with high enough density), some small amount of topology.

For the proof of the $C^1$ equivalence, we prove a decay property for density, coming
from a differential inequality like $\theta_x'(r) \geq C^{-1} r^{-1} [\theta_x(r)-\theta_x(0)]$.
Recall that the monotonicity of density (i.e., $\theta' \geq 0$) essentially follows by 
comparing $E$ with a cone.

For the differential inequality, we assume that $\theta_x(r)>\theta_x(0)$,
and construct a competitor (a deformation of $E$) which is significantly better than the cone 
over $E \cap \d B(x,r)$, and then hope that we get the right inequality.

One of the main ingredient is that for flat enough surfaces, $\Delta f = 0$ is a good 
approximation of the minimal surface equation, hence the graph of a harmonic extension 
often has a significantly smaller energy than the cone (the graph of the radial extension).

We use this, suitable gluing arguments between flat surfaces, and reduce to a length estimate 
(the ``full length property'' of $X$) on small perturbations of the geodesics of 
$K = \d B(0,1) \cap X$, where $X$ is the minimal cone that approximates $E$ in $B(0,1)$. 
When $n=3$ the three known minimal cones satisfy this; when $n \geq 4$ we don't know.

\section{Sliding almost minimal sets} \label{s8}

Return to the sliding Plateau problem and almost minimal sets.
We would like to extend the theorem of J. Taylor to this context, and the next 
accessible case seems to be the description of sliding almost minimal sets of dimension $2$, 
near a point of the sliding boundary (because otherwise we may use Theorem \ref{t7.1}).

We'll rapidly restrict to the case when the sliding boundary $\Gamma$ is a $C^1$
curve, or even a straight line, but we will say a few words about the case when $n=3$,
$\Gamma$ is surface (or even a plane) and $E$ has to stay on one side.

Again, the dream would be to be able to do the following. First, list all the minimal cones
relative to the subject, i.e., all the sliding minimal cones of dimension $2$, with a sliding
boundary which is a line or a plane (when $n=3$). And then, for each such cone $X$, 
prove an analogue of Theorem \ref{t7.1} at points $x\in E$ where $X$ is a blow-up
limit of $E$. If we do this we get a classification, up to $C^{1+a}$ diffeomorphisms,
of all the singularities of sliding almost minimal sets. But we'll see that this is a little optimistic

\subsection{Sliding minimal cones} \label{s8.1}

Recall Definitions \ref{d2.14} and \ref{d3.2} for sliding deformations and sliding 
almost minimal sets relative to a sliding boundary $\Gamma$.

Before we start making lists, observe that every plain minimal set is automatically
a sliding minimal set, regardless of its position relative to $\Gamma$. This is because
we add contraints in the definition of a deformation, so we have less competitors, so
it is easier to be sliding minimal. The reader should not be shocked by this in the context 
of soap films: it is often easy to introduce a (wet!) needle through a soap film, with very little
disturbance.

Conversely, a sliding minimal set that does not meet
$\Gamma$ is automatically a plain minimal set, because the condition \eqref{2.19} is void.

Here we give examples of sliding minimal cones of dimensions $1$ and $2$, and we'll only mention 
the new (non plain) ones. We start with the simplest problems. When $d=1$ and $\Gamma$ is a
line, $E$ can also be a half line with its endpoint on $\Gamma$ and perpendicular to $\Gamma$,
or a $V$-set composed of two half lines that leave from a same point $p\in \Gamma$,
that make the same angle with $L$, but in opposite directions because they also make an
angle at least $2\pi/3$ with each other. That last condition is needed, as before, 
because otherwise we may replace a piece of $V$ by a piece of $Y$ with the same basis.
The verifications are not hard, and neither is the case when $d=1$ and $\Gamma$ is a higher 
dimensional vector space. 

The next simple case is when $d=2$, $n=3$, $\Gamma$ is the horizontal plane, and $E$, as well
as its deformations, is required to contain $\Gamma$ and lie in the (closed) upper half space.
We'll call this \ub{Fang's case} because what follows comes from
\cite{FangHolder} and \cite{FangC1}. The sliding minimal cones are then $\Gamma$, the union of
$\Gamma$ with a half plane perpendicular to $\Gamma$, and the union of
$\Gamma$ with a half cone of type $\bY$ perpendicular to $\Gamma$. 

This setting is not so ridiculous. It comes from the following initial problem. We start from a
region $\Omega$ bounded by a smooth bounded surface $\Gamma$, so that $\Omega$
is on one side of $\Gamma$, and we look at sliding almost minimal sets $E \subset \ol \Omega$,
associated to the sliding boundary $\Gamma$, and with the additional constraint
that $E \supset \Gamma$. That is, sliding deformations as defined as in Section \ref{s2.7},
but we also required that $\varphi_t(x) \in \ol \Omega$ when $x\in E$.
Sliding almost minimal sets are defined in terms of these deformations. 

There is a Plateau problem associated to this. Take a closed set $E_0 \subset \ol\Omega$
that contains $\Gamma$, and try to find a sliding deformation $E$ of $E_0$ (that 
contains $\Gamma$), such that $\H^d(E)$ is minimal. 
If such a set exists, it is clearly sliding minimal. But this is also true when $E$ is 
a solution of the  Reifenberg homology problem associated to the boundary set $\Gamma$ 
as in Section~\ref{s2.6}, but with the additional constraint that $E \subset \ol\Omega$.

Now it turns out that any blow-up limit of a sliding almost minimal set in this setting is
one of the minimal cones of Fang's case (after a rotation, so that the tangent plane
to $\Gamma$ becomes the horizontal plane, and $\Omega$ lies above the plane).

And this is a case where the boundary regularity program mentioned above works fine:
it is possible to prove that if $E$ is a sliding almost minimal set with the constraints above, 
then near every point of $E \cap \Gamma$, the set $E$ is equivalent through a $C^{1+a}$ diffeomorphism to one of the minimal cones above. See \cite{FangHolder , FangC1};
these results also yield new existence results, but let us not elaborate.

We should observe that when we do not require $E$ to contain $\Gamma$ and lie on one side, 
the situation becomes complicated again, and probably more complicated than when $\Gamma$ 
is a curve, because $E$ may have complicated ways of being tangent to $\Gamma$.

\ms
From now on, we assume that $d=2$ and $\Gamma$ is a line. Often we'll be in $\R^3$,
but this is not always needed. There are a few \ub{new sliding minimal cones}, i.e., 
that are not plain minimal.

The simplest ones are the \ub{half planes} bounded by $\Gamma$; we'll call $\bH(\Gamma)$ the 
collection of these half planes. The verification that every $H \in \bH$ is minimal is not hard. 
Let $\pi$ denote the orthonormal projection on the plane $P$ that contains $H$, and let $F$ 
be any (sliding) deformation of $H$ in a large ball $B$ centered on $\Gamma$. 
With a little bit of topology (extend the deformation mapping to $P$ by symmetry
and use monodromy or degree theory), it is possible to show that $\pi(F\cap B)$ contains
$H \cap B$. Then $\H^2(F\cap B) \geq \H^2(\pi(F\cap B)) \geq \H^2(H\cap B)$, and the minimality
of $E$ follows. See Section 39 of \cite{C1W}. 

Next are the \ub{cones of type $\bV$}, or $\bV$-sets. 
Those are the unions $V = H_1 \cup H_2$ of two half planes $H_i \in \bH(\Gamma)$, 
that make an angle at least $2\pi/3$ along $\Gamma$.
This is the generalization of the $V$-sets of dimension $1$ above, the reason for the angle condition
is the same as above (otherwise, pinch), and the verification of minimality is not so hard either
(with a slicing argument; see \cite{CavaT} or Section 39 of \cite{C1W}). 
We'll call $\bV(\Gamma)$ the collection of these cones.

A special case of set $V \in \bV(\Gamma)$ is a plane that contains $\Gamma$.
It is also a special case of plain minimal cone, but we expect a different behavior from
$E$ near the points where such a plane is a blow-up. The same thing is (sadly) true of the cones 
$Y \in \bY(\Gamma)$, the cones of type $\bY$ whose singular set coincides with $\Gamma$.

In addition to these cones, Xangyu Liang suggested that the cone $\bQ$ 
over the edges of a cube, with great diagonal $\Gamma$, is also sliding minimal.
Soap experiments suggest that this may be true, but we have no proof in either direction.

There may be other sliding minimal cones that we did not think about, but the 
author would bet that (in $\R^3$) this is not the case.
In $\R^n$, $n \geq 4$, there are probably many 
more sliding minimal cones. Anyway, for the discussion below, $\R^3$ is enough trouble already.

\subsection{Regularity attempts near Sliding minimal cones} \label{s8.2}

Here $E$ will be a sliding minimal of dimension $2$ in $\R^n$, bounded by a smooth curve $\Gamma$.
In fact, we'll just assume that $\Gamma$ is a line, and say that the general case would be similar.
We want to follow the program mentioned above, and give a good local description of $E$ 
near any point $x\in E \cap \Gamma$, based on the knowledge of a blow-up limit $X$ of $E$ at $x$.
As one could guess, this will be harder to do for some cones $X$, and we shall not be able to complete
our program in all cases. The results below are taken from \cite{C1W}.

We start with the simplest case when $X$ is a \ub{half plane}. Then we get a full theorem of Taylor type.

\begin{mythm} \label{t8.1}
Let $E$ be a (coral) sliding minimal set in $B(0,2r_0)$, with a sliding boundary which is a line
through $0$. Suppose that $h(r) \leq \varepsilon (r/r_0)^{1+\alpha}$ for $0 < r < 2r_0$
and $d_{0,2r_0}(E,H) \leq \varepsilon$ for some half plane $H \in \bH(L)$.
Then $E$ is $C^{1+a}$-equivalent to $H$ in $B(0,r_0)$.
Here $\varepsilon >0$ and $a > 0$ are two small constants that depend only on 
$n$ and $\alpha > 0$.
\end{mythm}

Of course, by scale invariance, we could have restricted to $r_0 =1$.
Another way to put the $C^{1+a}$-equivalence would just have been to say that $E \cap B(0,r_0)$
is composed of a single $C^{1+a}$ face, bounded by $\Gamma$.
Anyway, $E$ has no holes, or complicated topology, in $B(0,r_0)$.

Notice that if $0 \in E \cap \Gamma$ and some blow-up limit of $E$
at $0$ is $H \in \bH(L)$, then the assumption of the theorem is satisfied for some radii $r_0$
that are as small as we want. Then the conclusion shows that in fact $E$ has a single blow-up limit
at $0$, and this limit is $H$. So that in fact the assumption is satisfied for all small $r_0$.

The ingredients of the proof are the following. 
If you are happy to settle with a bi-H\"older description (rather than $C^{1+a}$), 
a new monotonicity formula \cite{Mono} for balls $B(x,r)$ centered on $E \sm L$
(see \eqref{8.2} below), some compactness 
arguments using the almost-constant density principle relative to that monotonicity formula
(so that we can approximate $E$ by half planes or planes in many balls), the fact that $H$ is a cone of 
minimal density, and a Reifenberg parameterization if you want. 
This goes as in \cite{Holder} for a bi-H\"older variant of J. Taylor's theorem.

If you want the full $C^{1+a}$, you have to add a differential inequality, obtained by comparing
with minimal cones, some gluing argument to cut $E$ into two slightly smaller faces, and the fact 
that graphs of harmonic functions are usually closer to minimal than cones. The proof is not so different
than the one in \cite{C1}.

The monotonicity formula from \cite{Mono} concerns balls $B(x,r)$ that are centered 
on $E \sm \Gamma$.
For such balls, the usual density $\theta_x(r) = r^{-2} \H^2(E \cap B(x,r))$ is not always monotone
(for instance if $E \in \bH(\Gamma)$ and $x\in E \sm \Gamma$, because $\theta_x(r) = \pi$
for $r$ small and $\theta_x(r) \to \pi/2$ when $r \to +\infty$). We replace $\theta_x(r)$ with 
\begin{equation} \label{8.2}
F_x(r) = \theta_x(r) + r^{-2} \H^2(S_x \cap B(x,r)),
\end{equation}
where $S_x$ is the shade of $L$ seen from $x$, i.e. the half plane contained in the plane
that contains $x$ and $L$, and that lies behind $L$.

This new functional is nondecreasing, and is useful for some of the results of this subsection, 
because it is constant in some interesting cases (half planes or truncated $\bY$-sets). 
But not for full $\bY$-sets, which is the reason why these will cause trouble later.

\ms
The next simplest case is when $X$ is a \ub{generic cone $V$ of type $\bV$}.
Here generic means that the two half planes that compose $V$ make an angle
$\beta \in (\frac{2\pi}{3}, \pi)$. The statement is almost the same as above.

\begin{mythm} \label{t8.3}
Let $E$ be a (coral) sliding minimal set in $B(0,2r_0)$, with a sliding boundary which is a line
through $0$. Suppose that $h(r) \leq \varepsilon (r/r_0)^{1+\alpha}$ for $0 < r < 2r_0$
and $d_{0,2r_0}(E,V) \leq \varepsilon$ for some generic half plane $V \in \bV(\Gamma)$.
Then $E$ is $C^{1+a}$-equivalent to $V$ in $B(0,r_0)$.
Here $\varepsilon >0$ and $a > 0$ are two small constants that depend only on $n$, 
$\alpha > 0$, and $\beta \in (\frac{2\pi}{3}, \pi)$.
\end{mythm}

In particular, we need to make $d_{0,2r_0}(E,V)$ smaller when $\beta$ gets close to 
$\frac{2\pi}{3}$ or to $\pi$, because the sharp and flat cases below are a little more complicated.

Again, $E$ has no holes or complicated topology in $B(0,r_0)$; it is composed of two 
$C^{1+a}$ faces, both bounded by $\Gamma$, and that make at points $x\in \Gamma$ an angle
$\beta(x)$ that may vary slowly with $x$, but stays close to $\beta$.

The proof of this, and also of the next results, is again similar to what was done in \cite{C1}, 
but with the adapted monotonicity formula from \cite{Mono}. 
There is no simpler H\"older argument that works here, because part of the proof
consists in proving that the angle $\beta(x)$ that the two faces of $E$ make at 
$x\in \Gamma \cap B(0,r_0)$ varies slowly enough to avoid the more complicates minimal cones
studied below. 

\ms
The next case is when $X$ is a \ub{sharp cone of type $\bV$}, which means that
the two half planes that compose $V=X$ make an angle $\beta = \frac{2\pi}{3}$.
This is the first case where although we have a good description of $E$ in the 
$C^{1+a}$ category, we cannot say that $E$ is equivalent to $X$, and it even has a different
topology in general. We first give an incomplete statement, and then add some information.

\begin{mythm} \label{t8.4}
Let $E$ be a (coral) sliding minimal set in $B(0,2r_0)$, with a sliding boundary which is a line
through $0$. Suppose that $h(r) \leq \varepsilon (r/r_0)^{1+\alpha}$ for $0 < r < 2r_0$
and $d_{0,2r_0}(E,V) \leq \varepsilon$ for some sharp $\bV$-set $V \in \bV(\Gamma)$.
Then in $B(0,r_0)$, $E$ is composed of two main faces $F_1$ and $F_2$, plus maybe a third
thin set $F_3$ composed of one or more ``vertical faces'', that meet along a curve $\gamma$
that is in general partially contained in $\Gamma$.
Here $\varepsilon >0$ and $a > 0$ depend only on $n$ and $\alpha > 0$.
The two main faces meet all along $\gamma$, and $F_3$ is bounded on one side
by $\gamma \sm \Gamma$, where it meets $F_1$ and $F_2$ with $\frac{2\pi}{3}$ angles,
and on the other side by $\Gamma \sm \gamma$, where it is attached to $\Gamma$ as in
Theorem~\ref{t8.1}.
\end{mythm}

\begin{figure}[!h]  
\centering
\includegraphics[width=10.cm]{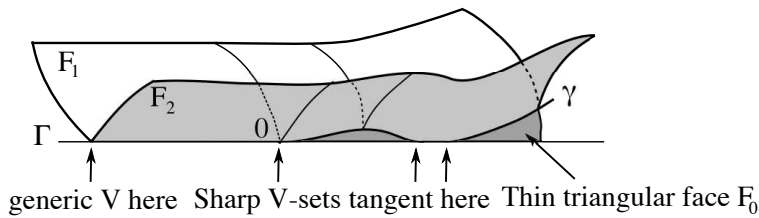}
\caption{$E$ near a sharp cone of type $\bV$ \label{fV}}
\end{figure}

See Figure \ref{fV}. 
The set $E$ looks like a $\bV$-set, but which we can pinch a little in some places along $\Gamma$
and also open slightly in some other places.
The general idea is that the two face $F_1$ and $F_2$ can escape from $\Gamma$,
but have to leave a thin wall $F_3$ that connects them to $\Gamma$. Along $\gamma \sm \Gamma$,
$E$ has a singularity of type $\bY$. Along $\Gamma \cap \gamma$, $E$ is described by
Theorem~\ref{t8.3}, with anges $\beta(x) \geq \frac{2\pi}{3}$ that may vary.
Finally, $\gamma$ is tangent to $\Gamma$ when it leaves it.
We do not exclude the possibility that $\gamma$ may leave $\Gamma$ and return to it
infinitely many times. Yet, for a minimal set, we would not expect anything like this, but just
that $F_1$, $F_2$, and $\gamma$ leave $\Gamma$ frankly.

It seems to me that Kenneth Brakke \cite{Br2} thought that this behavior may 
happen even when $V$ is a plane that contains $L$, but I claim that this is not the case.

\ms
Next we do a small digression. Suppose $0 \in E \cap \Gamma$, and that the blow-up limit
$X$ of $E$ at $0$ is one of the plain minimal cones (a set of type $\bP$, $\bY$, or $\bT$),
and that $X \cap \Gamma = \{ 0 \}$.
In this case what really happens is that the sliding condition does not really mean much, 
because in fact the set $E$ is transverse to $\Gamma$. Some verification needs to be done, but
nothing too hard. Then the proof of regularity for plain almost minimal sets applies in this case.
This is not shocking: you may take a soap film, and then quietly put a (wet, this is important!)
needle through it, and the set will essentially not be deformed; the needle will just cross $E$
transversally. The tangential case, discussed soon,  is a little more delicate, as the film
may prefer to follow the needle a tiny bit.

\begin{figure}[!h]  
\centering
\includegraphics[width=8.cm]{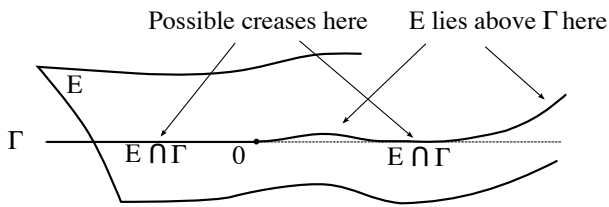}
\caption{$E$ near a plane that contains $\Gamma$. \label{fP}}
\end{figure}

\ms
We are ready for the next case, when $X$ is \ub{a plane $P$ that contains} $\Gamma$. 
Keep the same assumptions as in the theorems above, but with $d_{0,2r_0}(E,P) \leq \varepsilon$
for such a plane $P$. One possibility is just that $E$ is a smooth surface which is tangent
to $E$ at $0$, but a slightly more general behavior is allowed too, where $E$ is attached
to $\Gamma$ along an open set $I \subset \Gamma$, as in Theorem~\ref{t8.3}. 
See Figure \ref{fP}.
That is, $E \cap B(0,r_0)$ is still a Lipschitz graph over (its projection on) $P$,
smooth away from $\Gamma$, and that may have a small crease along $I \subset \Gamma$.
This happens with soap films (put a needle or a curve tangentially along the film,
and you should see some attraction), but an additional reason for this is probably capillarity,
which is not the subject here.

\begin{figure}[!h]  
\centering
\includegraphics[width=12.cm]{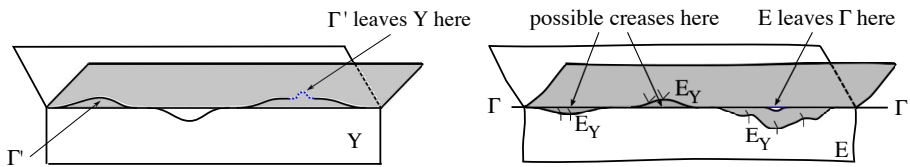}
\caption{$E$ near cone of $\bY(\Gamma)$. On the left, the image by $\psi^{-1}$, where $E$ becomes $Y$.
On the right, $E = \psi(Y)$ and $\Gamma = \psi(\Gamma')$. 
We tried to mark directions along the $E_Y$ of $\bY$-points of $E$.
\label{fYY}}
\end{figure}

\ms
Let us skip some other cases that are similar to the previous ones (see \cite{C1W}), and go directly to
\ub{the main bad case}, when $X$ is \ub{a cone $Y \in \bY(\Gamma)$}, i.e., which contains $\Gamma$.
There is a natural conjecture about this case, which is that $E$ is the image
of $Y$ by a homeomorphism $\psi$, which sends some curve 
$\Gamma' = \psi^{-1}(\Gamma)$ to $\Gamma$, but with no special requirement 
concerning the position of $\Gamma'$ towards $Y$. The way $E$ organizes itself around $\Gamma$
(and in particular the parts of $\Gamma$ where $E$ looks like a $\bY$-set, a $\bV$-set, a $\bH$-set,
or $\emptyset$) would then follow from the position of $\Gamma'$ relative to $Y$, as in the the 
previous case when $X \in \bV(\Gamma)$. See Figure \ref{fYY}.

Unfortunately, this seems harder to prove than the previous results. The main problem is that the
monotonicity formula from \cite{Mono} does not give good results in the case of $\bY$-sets
(it is strictly increasing in the case of interests), and it is not clear that there is a replacement for it.
Another way to state the issue is to say that since we expect singularities of type $\bY$ in
$E \sm \Gamma$, having a good monotonicity formula for balls centered at these points should
help a lot, especially if we intend to apply a form of the almost-constant density principle 
to control the geometry of $E$ in balls centered at these points.

We end this section with a maybe not so easy exercise.
We know one way in which the sliding almost minimal set $E$ associated to the smooth 
boundary curve $\Gamma$ may leave $\Gamma$ (i.e., at a point where $E$ las a tangent plane 
that contains the tangent line to $\Gamma$). Are the other ways to do this?

\subsection{Further questions} \label{s8.3}

For me the next main question is the regularity of the sliding almost minimal set $E$
near a point where $E$ admits a blow-up limit which is a $\bY$-set $Y$ that contains $\Gamma$.
If we can do this, there is a good hope that a general regularity result (maybe only in the H\"older
category, near some isolated points of $E \cap \Gamma$) will ensue, and then also that existence
results will follow.

Making sure that we have a full list of sliding minimal cones in $\R^3$ (when $\Gamma$ is a line)
would be nice too. And if this is too easy, do the same thing in $\R^4$.

And then we may also worry about the following. Let $\Gamma \subset \R^3$ be a smooth 
closed curve, and choose a closed set $E_0$, if possible well attached to $\Gamma$ so that
the sliding Plateau problems below are nontrivial.
For $\varepsilon > 0$ small, let $T_\varepsilon$ be a smooth open tube around $\Gamma$, 
of width roughly $\varepsilon$. 
Consider Fang's problem where we minimizer $\H^d(F)$ among sliding deformations $F$ of 
$E_\varepsilon = \d T_\varepsilon \cup E_0 \sm (T_\varepsilon)$,
inside $\Omega = \R^3 \sm T_\varepsilon$, and with the boundary $\d T_\varepsilon$.
 
How do the solutions of this problem (which exist by \cite{FangC1}) 
converge to solutions of the sliding Plateau problem associated to $\Gamma$ and $E_0$,
assuming that there are some?
The reader should not get confused like the author: the extra condition that $E \supset \Gamma$
does not mean much when $\Gamma$ is one-dimensional, and to some extent it should have less
and less importance for the Fang problem on $\R^3 \sm \d T_\varepsilon$, except to give 
the existence of minimizers.

Probably it is better to solve the problem about the bad $\bY$-cone, and then the existence
problem for the sliding Plateau problem, before we do this.

We did not talk at all about what happens to the regularity results when we replace $\H^d$
with other elliptic integrands, possibly not even smooth in one of the variables. A simple example
of this (but for which only partial descriptions exist at this time) is to use the functional
\begin{equation} \label{8.5}
J(E) = \H^d(E \sm \Gamma) + \alpha \H^d(E \cap \Gamma)
\end{equation}
for some $\alpha \in [0,1)$, for instance in the context where $\Gamma$ 
is the horizontal plane in $\R^3$ and $E$ is a sliding almost minimal set with boundary 
$\Gamma$, with the additional constraint that $E$ lies un the (closed) upper half space. 
See \cite{CavaT} for some initial partial results in this context, and \cite{FangEx, Ita3} 
and probably many others in the general context.

\section{V. Feuvrier and existence results}

In this last section, we present a general scheme for proving existence results.
The main tool here will be a construction \cite{FeuPoly} of adapted polyhedral networks, and then some
of the results mentioned above.

This is one of the three or four systematic ways to get existence results that the author knows of;
maybe not necessarily the simplest (this also depends on people's personal knowledge), but
the author thinks it is quite nice and natural, and was perhaps misunderstood so far.
It is somewhat close to Reifenberg's initial approach in \cite{Reifenberg}, or the proof of
\cite{HaPu}, but there are significant differences.

Another type of proof is through currents (or flat chains, or varifolds) and the celebrated
compactness theorem; traditionally they tend not to work so well for sets and size minimizing
currents, but there is an approach by Almgren (with varifolds and flat chains, for the Reifenberg 
homology problem with integrands), which was recently nicely put together by U. Menne,
Y. Fang, and S. S. Kolas\'{i}nski \cite{FaKo}. 

And finally there is a beautiful recent approach by De Lellis, Ghiraldin, and Maggi \cite{Ita1}
and De Philippis, De Rosa, and Ghiraldin \cite{Ita2}, where they take a minimizing sequence
of sets $E_k$, consider a weak limit of the measures $\mu_k = \H^d_{\vert E_k}$,
and use it to find a good minimal set. 

It is interesting to notice that the different approaches have similar or parallel points, but also
differences (so that they don't always work in the same circumstances).

\subsection{Presentation} \label{s9.1}

The method that we describe here consists in taking a suitably improved minimizing 
sequence of sets, taking a subsequence that converges for the Hausdorff distance, 
and hopefully the limit will be the desired minimizer. 
It was used by V. Feuvrier \cite{FeuThesis, Feu3} (to get sliding $2$-dimensional
minimal sets in a situation without boundary), X. Liang \cite{LiangTopo} (to get minimal sets 
under topological constraints \cite{LiangTopo}, and Y. Fang \cite{FangEx} (to get Reifenberg
homological minimal sets under general conditions). The first general presentation (too vague?)
dates from \cite{Purdue}, which apparently was not convincing enough. Yet the author does not
intend to give all the details. 

The advantage of Feuvrier's method is that it looks simple, and also some 
of the ideas are pleasant. But there are technical problems in some of the ingredients
of this proof, more or less as with the other proofs.

We'll do things in the setting of \cite{FeuThesis, Feu3}, because it is the simplest, and in particular
we will not get distracted by boundary problems or topology. 
Let us work on a flat compact connected manifold $\cM$, without boundary, of the 
following type: $\cM$ is obtained from a finite number of dyadic
cubes, by identifying some faces. The point of this setting (which indeed looks like a toy model)
is the following. First, the absence of boundary will allow us to forget about any problem of sliding
boundaries on other ways to state Plateau problems. These would add technical annoyances,
but not necessarily bad ones. 

Next, the simple flat structure will allow us to use dyadic grids on $\cM$. With a regular manifolds, 
we would need to worry about how to define polyhedral grids in an environment where the change 
of chart mappings are not necessarily affine. Apparently this technical difficulty is being taken 
care of by Feuvrier (constructing adapted grids on a manifold), but this adds to the 
technical complication of the construction.

Also, we'll work in dimension $d=2$; this will allow us to use J. Taylor's regularity theorem
or its extension in \cite{Holder} in the last step of the verification. But we'll pretend to work
in dimension $d$ for some time, up to the moment when we have to take $d=2$.

We give ourselves a closed set $E_0 \subset \cM$, with $\H^d(E_0) < +\infty$, and
consider the class $\cE$ of images of $E_0$ by continuous deformations. That is, $E \in \cE$
if there is a (continuous) one parameter family $\{ \varphi_t \}$ of continuous mappings
$\varphi_t : E_0 \to \cM$, with $\varphi_0(x) = x$ on $E_0$, and $E = \varphi_1(E)$.
And we want to minimize $\H^d$ in the class $\cF$, i.e., find $F \in \cE$ such that
\begin{equation} \label{9.1}
\H^d(F) = m, \ \text{ where } m = \inf_{E \in \cE} \H^d(E).
\end{equation}

Of course this is more interesting if $m > 0$, which excludes the case when $E_0$
can be contracted in $\cM$ to a point (or a lower dimensional object). This is why we want to 
take a manifold $\cM$ with some topology. 
The problem could turn out to be a little too easy to solve
with our flat structure, but replacing $\H^d(E)$ with $\int_E f(x) d\H^d(x)$ for some 
H\"older-continuous function $f$ on $\cM$ such that $C^{-1} \leq f \leq C$, or even
one of the slightly more general integrands of \eqref{3.4}-\eqref{3.5},
could make the problem less trivial without really changing the proof.

We describe a general scheme now, which can work in other circumstances as well.
Since we don't know any other way, we'll take a minimizing sequence $\{ E_k \}$, i.e., such that
\begin{equation} \label{9.2}
E_k \in \cE \ \text{ and } \, \lim_{k \to +\infty} \H^d(E_k) = m.
\end{equation}
The logical thing to do next is try to get a subsequence that converges to a limit $F$,
and then try $F$ in \eqref{9.1}. It would be nice if we could use a strong (fine) topology for this,
because this would probably make it simpler to prove the minimality of $F$. For instance,
it would be great if we could write $E_k = \varphi^{k}_1(E_0)$ for some
sequence of one parameter mappings $\{ \varphi^k_t \}$, $0 \leq t \leq 1$,
and on top of this the  $\varphi^k_t$ converge uniformly on $E_0 \times [0,1]$,
but we should not dream too much: how are we going to get enough compactness for
such a sequence? Notice however that this is a little bit what Rad\'o and Douglas did
(see Section \ref{s2.2}), by restricting to smaller classes of parameterizations, but then there are 
other issues and we no longer want to do this.

Anyway, here we decide that we want to make sure that convergent subsequences exist,
so we take the weakest topology on $\cE$, the topology of local Hausdorff convergence 
that was defined in near \eqref{5.3}. But here $\cM$ is compact, so this is just the topology defined
by the usual Hausdorff distance.

\subsection{The need for a quasiminimal haircut} \label{s9.2}

There are two \ub{obvious problems} with what we just started. The first one is \ub{hairs}.
If we take a limit of a subsequence of $\{ E_k \}$, it could converge to very large sets, 
even with infinite measure, that are certainly not minimizers. Indeed, starting with two-dimensional 
sets in $3$-space, it is often very easy to replace any good-looking set $E_k$ with an uglier set 
$\wt E_k$, which is a deformation of $E_k$ obtained by growing long thin hair from it, in such a way
that $\wt E_k$ is $\varepsilon$-dense in $\cM$ for any small $\varepsilon > 0$ given in advance,
and yet $\H^d(\wt E_k) \leq \H^d(E_k) + \varepsilon$. If we have the bad fortune to use $\wt E_k$
instead of $E_k$, it will converge to the whole $\cM$, something that we want to avoid.

So we want to make $E_k$ cleaner before we take a Hausdorff limit. This is what we shall call
doing a haircut; in his time Reifenberg did this too, but maybe we want to be more systematic.
We'll find a new minimizing sequence $\{ F_k \}$ in $\cE$, obtained from $\{ E_k \}$,
and then only extract a subsequence that converges to a limit $F_\infty$.

Our second, in fact related problem, is the \ub{convergence of Hausdorff measure}. 
We want make sure that for this new sequence $\{ F_k \}$ that converges to $F_\infty$,
\begin{equation} \label{9.3}
\H^d(F_\infty) \leq \liminf_{k \to +\infty} \H^d(F_k) = 
\lim_{k \to +\infty} \H^d(F_k) = m.
\end{equation}
And in fact the main point of the haircut is to have \eqref{9.3}.
Now we have seen such a property before, in Lemma \ref{t5.9},
when we were studying limits of reduced almost minimal sets. 
Now it will be a little too hard to find a sequence $\{ F_k \}$, as above, which is 
composed of almost minimal sets; in fact the best way to get almost minimal sets
seems to be to minimize a variant of $\H^d$, and this is precisely what we were trying 
to do in the first place.
But fortunately, the notion of quasiminimal sets (defined soon) is just right: it will be easier
to produce a sequence of quasiminimal sets, and yet the notion is sufficiently strong for
\eqref{9.3} to hold as soon as the sets $F_k$ are reduced quasiminimal sets, 
with quasiminimality constants that does not depend on $k$.
So let us define quasiminimal sets. We keep the sliding boundary condition to make
a more general condition, but in the setting that we chose, there is no sliding boundary
and we can take $\Gamma = \emptyset$. We also localize to an open set $U$
for the same reason.

\begin{mydef} \label{d9.4}
Let $U \subset \R^n$ be open and let $\Gamma \subset \R^n$ be a closed set.
We say that $E$ is a (sliding) quasiminimal set in $U$, with the sliding boundary 
$\Gamma$, at the size $\delta_0$ and with the quasiminimality constant $M$, 
when for each closed ball $B \subset U$ of radius at most $\delta_0$ and each 
deformation $\{ \varphi_t \}$, $0 \leq t \leq 1$, for $E$ in $B$
(as in Definition \ref{d2.14}), we have the following estimate. Set 
\begin{equation} \label{9.5}
W = \big\{ x\in E \, ; \, \varphi_1(x) \neq x \big \} \subset E \cap B;
\end{equation}
then 
\begin{equation} \label{9.6}
\H^d(E \cap W) \leq M \H^d(\varphi_1(E \cap W)).
\end{equation}
\end{mydef}

\ms
The definition is a little complicated, but makes sense.
It says that we may make $E$ smaller by deforming a part of $E$, but we are not able 
to save more than $99\%$ of what we modify, say. 
And it makes sense to allow deformations that modify just a tiny bit of $E$ in a ball. 
The definition is essentially Almgren's \cite{AlMemoir}
(he called these sets ``restricted sets''), but in \cite{DSPacific} and \cite{DSMemoir}
we found out that we wanted exactly that definition.

In the plain case, we don't bother with the existence of a whole one parameter family 
$\{ \varphi_t \}$, because they would be easily reconstructed from $\varphi_1$, 
and we state the definition directly in terms of $\varphi_1$.
Also, we could allow other compact sets than balls as the sets where the deformation
takes place (this could possibly make a difference on complicated domains $U$ where 
large closed balls are difficult to find, but so far the author did not see the difference.
Finally in \cite{Sliding} the author felt compelled to add one more type of objects, 
the generalized quasiminimal sets, where one replaces \eqref{9.6} by 
\begin{equation} \label{9.7}
\H^d(E \cap W) \leq M \H^d(\varphi_1(E \cap W)) + r^d h(r),
\end{equation}
where $r$ is the radius of $B$ and $h$ is a gauge function as above.
This gives a little more generality, and the same proofs apply anyway.

As was hinted before, the results of Sections \ref{s4} and \ref{s5} are still true with
quasiminimal sets, with essentially the same proofs. This includes the results about limits,
but not the more precise results that involve monotonicity or even less epiperimetric 
inequalities. This is the main point of \cite{Sliding}.

Now how much more flexible is the theory of quasiminimal sets? 
In dimension $d=1$, (connected) quasiminimal sets are essentially the same
as Chord-arc curves (length of the arc $\Gamma(a,b)$ less than $C|b-a|$), which,
when $n=2$ and for unbounded curves, are the bilipschitz images of a line.

It is easy to see that (in the plain case and, say, in $\R^n$ to avoid
complications near the boundary) the image of a quasiminimal set by a bilipschitz 
mapping is also quasiminimal, maybe with a larger constant $M$ and a smaller $\delta_0$.
This includes bilipschitz images of (locally) minimal sets, which happen to be 
(locally) quasiminimal with constant $M=1$; the verification in unfortunately 
less pleasant as it should be, due to the strange definition with the set $W$, but the reader
should not be surprised either.

Anyway, Lipschitz graphs, for instance, are quasiminimal, and this gives a vague idea of
their regularity, since it is known that (away from the sliding boundary) they are 
uniformly rectifiable.

Another way to get a quasiminimal set is to minimize a functional $J(E) = \int_E f(x) \d\H^d(x)$
(instead of $\H^d(E)$), where we only know that $C^{-1} \leq f(x) \leq C$ (and in particular
we have no good continuity property for $f$). We could even let $f$ be an integrand that
depends on the approximate tangent plane to $E$ (and is defined some other way on the unrectifiable part of $E$). The verification is easy.

This may look a little bad, since then we cannot expect quasiminimal sets to have 
better regularity properties than Lipschitz graphs, but this is also what will give us enough 
flexibility for the proof of existence that we try to describe here.

So let us return to our existence problem. We decide that we'll do a quasiminimal haircut,
i.e., find a new minimizing sequence $\{ F_k \}$, such that
\begin{equation} \label{9.8}
\begin{aligned}
&\text{the sets $F_k$ are quasiminimal, with some constants}
\cr
&\hskip.8cm \text{$\delta_0 >0$ and $M \geq 1$ that do not depend on $k$.} 
\end{aligned}
\end{equation}
Let $F_k^\ast$ denote the core of $F_k$; our life is simpler if the $F_k$ are reduced
(i.e., $F_k^\ast = F_k$, but in fact we do not really expect this, and we'll have to discuss
the difference). Anyway, a small verification (a little less pleasant in the sliding case, but
this fact is checked in \cite{Sliding}) shows that the $F_k^\ast$ are quasiminimal,
with the same constants $\delta_0$ and $M$. Let us replace $\{ F_k \}$ by a subsequence 
for which the $F_k^\ast$ converge (locally for the Hausdorff distance, as usual) to a limit
$F_\infty^\ast$. The theorem about limits says that
since \eqref{9.8} holds for the $F_k^\ast$ and at the same time $\{ F_k \}$ is a minimizing
sequence (and hence the $F_k^\ast$ cannot be improves by deformation by more than 
numbers $\varepsilon_k$ that tend to $0$; again this takes a small amount of checking,
to see that the potential deformation can be extended to $F_k \sm F_k^\ast$, at no cost),
we get that 
\begin{equation} \label{9.9}
F_\infty^\ast \ \text{ is a minimal set} 
\end{equation}
(if we were working on a different problem, with a boundary condition, 
it would be sliding minimal), and, as in Lemma \ref{t5.9},
\begin{equation} \label{9.10}
\H^d(F_\infty^\ast) \leq \lim_{k \to +\infty} \H^d(F_k^\ast)
= \lim_{k \to +\infty} \H^d(F_k) = m
\end{equation}
(because $\H^d(F_k \sm F_k^\ast) = 0$ and by \eqref{9.3}). 
We claim that our life will be a little simpler then, but how do we do our haircut?

\subsection{dyadic grids, polyhedral nets, and quasiminimal haircuts.} \label{s9.3}

We start with a (slightly too) simple idea: we start from our initial set $E_k$
(from the minimizing sequence), and we perform a Federer-Fleming projection
on a dyadic grid of small mesh size. We get a new set $E_k^1$, which is contained in a 
finite union of faces of dimension $d$. Some faces $S$ are not entirely contained in $E_k^1$,
and when this happens, we can project again on the boundary of $S$, getting a new set
$E_k^2$ composed of full $d$-faces and subsets of $(d-1)$-faces. This is in fact enough
for us, but for the sake of organization, let us continue all the way to the set $E_k^d$,
which is still a deformation of $E_k$ (so it lies in $\cE$), and in addition lies in the class
$\cG_k$ of finite unions of faces of dimensions at most $d$ of our dyadic grid. 

We like finite problems, because we can solve them. 
So we replace $E_k^d$ by a new set $F_k^0$, which we choose to minimize $\H^d(F_k)$ 
in the class $\cE \cap \cG_k$. This may be quite a brutal change, so what was the point 
of the $E_k^d$?
We hope to say that $\{ E_k^d \}$ too is a minimizing sequence, i.e.,
\begin{equation} \label{9.11}
\lim_{k \to +\infty} \H^d(E_k^d) = m
\end{equation}
and then it will follow that $\{ F_k \}$ is minimizing too, since
\begin{equation} \label{9.12}
\H^d(F_k) \leq \H^d(E_k^d).
\end{equation}
We will return soon  to the issue of whether \eqref{9.11} holds or not, but first 
let us say what we win by replacing $E_k^d$ by $F_k$. The point is that \eqref{9.8} holds,
just because $F_k$ minimize $\H^d(F_k)$ in the class $\cE \cap \cG_k$.
We can take $\delta_0 = 1$ (in fact, provided that earlier, we took a dyadic net of
mesh size at most $n^{-1/2}$, say), and the constant $M$ does not depend on 
$k$, just on $n$ and the fact we use (small enough) dyadic grids.
And of course it is important that estimates do not depend on $k$ and do not get worse 
when we take smaller dyadic grids.

Let us see how the proof starts. 
Let $G = \varphi_1(F_k)$ be a competitor for $F_k$, as in Definitions \ref{d9.4}
and \ref{d2.14}. Of course $G$ does not necessarily lie in the class $\cG_k$, but we can
apply a Federer-Fleming projection to it, as we did earlier to get $E_k^d$, to get a competitor
$G'$ for $G$ that lies in $\cG_k$. By definition, $\H^d(F_k) \leq \H^d(G')$.
Then we look carefully at how the points move, notice that $G'$ is almost as good a
competitor as $G$ was (because whenever we move points from the interior of a face
to the boundary of that face, we multiply its measure by at most $C$), do the algebra, and
conclude. See Chapter 11 of \cite{DSMemoir}, and in particular Proposition 11.13
on page 90; in the present situation, $F_k$ is even a restricted minimizer.

We return to \eqref{9.11}. Unfortunately, things are not so simple; even for $d=1$
in $2$-space, $E_k$ could coincide with the first diagonal in a large ball, and when we
project it on a dyadic grid (with faces parallel to the axes), we multiply its length by
roughly $\sqrt2$, which of course is bad and ruins our chances for \eqref{9.11}.
This is where the polyhedral nets constructed by V. Feuvrier in \cite{FeuPoly}, 
and that we shall describe now, will be useful.

Here our assumption that our space $\cM$ is simple (essentially, a Euclidean space with
identifications) will simplify our life. Suppose we are given a finite collection of disjoint sets
$U_i \subset \cM$, such that $\dist(U_i, U_j) \geq \varepsilon$ for $i \neq j$ and some small
constant $\varepsilon > 0$. Then suppose that we chose for each $i$ a dyadic grid,
whose mesh size and direction are allowed to depend on $i$. What Feuvrier provides is a
decomposition of $\cM$ into small convex polyhedra, with the usual structure of dyadic cubes
concerning the decomposition into faces, and whose restriction to an 
$\varepsilon/10$-neighborhood of $U_i$ coincides with a dyadic subnet of the initial net
that was given. The size of the polyhedra may be very small, but what is important is that
the various angles between faces are bounded from below, by a positive number that depends
only on $n$ (Feuvrier says that the faces are uniformly ``rotund''). Because of this,
we can perform Federer-Fleming projections on these polyhedral nets, and the various 
constants in the constructions (coming from the constants in \eqref{4.9}, \eqref{4.10},
and the like) are uniform.

The construction of these nets is somewhat painful, but the result is not too shocking,
even though this is probably the less amusing part of the proof described here.

Now the reader guessed what we shall do. For each $k$, we want to use a Feuvrier net adapted
to our original candidate $E_k$. Let us assume that $E_k$ is rectifiable (otherwise,
we could also project on dyadic cubes in some balls, to make almost all of its unrectifiable 
part disappear, and modify the rest of the argument below). Then almost-cover $E_k$
(using the Vitali covering argument as for Lemma \ref{t5.9}) by disjoint small balls $B_j$ where
$E_k$ is very well approximated by $d$-planes $P_j$. For each of these balls, take a dyadic 
grid such that $P_j$ is a coordinate $d$-plane, so that when we do a Federer-Fleming
projection of $E_k \cap B_j$ in this dyadic grid, we add almost no measure. Complete
the grid as in \cite{FeuPoly}; in $\cM \sm \cup_j B_j$, we may multiply the measure
by $C$ when we project on the $d$-faces, but we don't care, because the total mass
of $E_k$ there is as small as we want. So, with this adapted grid, we manage to do our first
Federer-Fleming projection, replace $E_k$ with $E_k^d$ as we did before, and in addition
make sure that $\H^d(E_k^d) \leq \H^d(E_k) + 2^{-k}$. That is, we get \eqref{9.11}.

The grid has a mesh size that will probably tend to $0$ when $k$ tends to $+\infty$,
but we don't care. What is important is that the rotundity, then the Federer-Fleming constants,
and finally the quasiminimality constants for $F_k$ in \eqref{9.8}, still obtained by the proof
of Proposition 11.13 11 of \cite{DSMemoir} (for instance), do not depend on $k$.

At this point we realized our quasiminimal haircut (i.e., replaced $\{ E_k \}$ with a
new sequence $\{ F_k \}$ which is still minimizing (by \eqref{9.11} and \eqref{9.12}),
and satisfies the quasiminimal condition \eqref{9.8}. As was said before, we can replace
$\{ F_k \}$ by a subsequence for which the cores $F_k^\ast$ converge to a limit
$F_\infty^\ast$, and we get the two properties \eqref{9.9} and \eqref{9.10}.

\subsection{End of game.} \label{s9.3}

So far we used very little (the most unpleasant part is our assumption on $\cM$,
which probably can be removed as soon as Feuvrier constructs his cubes in a manifold).
We could even work with elliptic integrands (because Lemma \ref{t5.9} works in this context).
But it is important to notice that we may not be finished yet.

Indeed, although \eqref{9.9} says that $F_\infty^\ast$ is a minimal set 
(and \eqref{9.10} that it has the right measure), we do not know yet whether it lies
in the class $\cE$, or the class where we wanted to minimize $\H^d$. 
This is where we need to adapt our proof with the specific problem at hand.

In general we want to minimize $\H^d(E)$, or a variant, in a class $\cE$ such that
deformations of $E \in \cE$ automatically lie in $\cE$. This way, we start with sets 
$F_k \in \cE$. In some cases this may also imply that the core $F_k^\ast$ automatically
lies in $\cE$, but not always. And then we take a limit of the sets $F_k^\ast$, and once more
it could be that $\cE$ is not stable under limits.

Let us first say how to deal with these two issues in the setting of that was described
at the beginning of this section, and then, we'll rapidly discuss other cases.

We decided to take an initial set $E_0$, and for $\cE$ the class of (continuous) deformations
of $E_0$ in our test manifold $\cM$, and also to restrict to $d=2$. The obvious difficulty
is that since we are in fact dealing with parameterizations (by $E_0$), taking limits should
be hard, because we have no control on the regularity of the parameterizations.
What will save us is the regularity of the limit.

By \eqref{9.9}, $F_\infty^\ast$ is a minimal set. We also organized our problem so that there
is no boudary $\Gamma$ to complicate matters. We start with the simpler case 
when $\cM$ is of dimension $3$. 
Then we can apply Theorem \ref{t7.1}, and get that locally $E$ is a $C^{1+\alpha}$ 
variant of a minimal cone of type $\bP$, $\bY$, or $\bT$. 
It is also compact (because $\cM$ is), so we can use this regularity to build a 
Lipschitz retraction near $\cM$, which we can even glue to the identity
(far from $E$), to get a (continuous) one-parameter family of continuous mappings 
$\psi_t$, $0 \leq t \leq 1$, from $\cM$ to $\cM$, with the following properties. 
First, there is an $\varepsilon > 0$ such that
\begin{equation} \label{9.12}
\psi_t(x) = x \ \text{ for $0 \leq t \leq 1$ when } \dist(x,F_\infty^\ast) \geq 2\varepsilon, 
\end{equation}
\begin{equation} \label{9.13}
\psi_1(x) \in F_\infty^\ast \ \text{ when } \dist(x,F_\infty^\ast) \leq \varepsilon,
\end{equation}
and (coming from the $C^{1+\alpha}$ regularity,
\begin{equation} \label{9.14}
\psi_1 \ \text{ is Lipschitz.} 
\end{equation}
The reader imagines also that $\psi_t(x) = x$ for $x\in F_\infty^\ast$; we 
shall not really need that, but we get it from the construction anyway.
The proof is not hard; the point is that the existence of a Lipschitz retraction turns
out to be local (we can compose retractions).

Now we are ready to construct the desired minimizer for $\H^2$ in $\cE$. 
Recall that $F_k^\ast$ tends to $F_\infty^\ast$, hence for $k$ large, $F_k^\ast$ 
is contained in an $\varepsilon$-neighborhood of $F_\infty^\ast$. We claim that 
$F = \psi_1(F_k)$ does the job. First of all, $F \in \cE$ because 
$F$ is a deformation of $F_k$, which itself lies in the class $\cE$ by construction.
Since $\psi_1$ is Lipschitz, $\H^d(\psi_1(F_k \sm F_k^\ast)) = 0$ because
$F_k \sm F_k^\ast$ is a finite union of faces of dimensions at most $d-1$.
So $\H^d(F) = \H^d(\psi_1(F_k^\ast) \leq \H^d(F_\infty^\ast) \leq m$, by 
\eqref{9.13} and \eqref{9.10}. So $F$ is the desired minimizer.
In fact $\psi_1(F_k^\ast) = F_\infty^\ast$, because
otherwise the inequality above would be strict, so the core of $F$ is our minimal
set $F_k^\ast$, but we did not really need to know this.

The argument also works when $\cM$ is of dimension larger than $2$, with minor differences.
The main one is that we do not know the exact list of minimal cones of dimension $2$
in $\R^n$, and \cite{C1} only gives a local biH\"older description of $F_\infty^\ast$,
with any H\"older exponent smaller than $1$. Yet, the combinatorial description of
minimal cones (and then of $F_\infty^\ast$) in terms of faces that meet with $120^\circ$
angles is still enough to build a retraction near $F_\infty^\ast$, and then mappings
$\psi_t$ as above, except that \eqref{9.14} no longer holds, and has to be replaced
by a H\"older condition, fortunately with an exponent which is as close to $1$ as we want.
Then we can complete the argument as above, arguing now that 
$\H^d(\psi_1(F_k \sm F_k^\ast)) = 0$ because $\psi_1$ is H\"older and 
$F_k \sm F_k^\ast$ is at most $(d-1)$-dimensional.
And in this context, we make sure not to require that the final mappings of deformations
are Lipschitz in the definition of $\cE$, because otherwise we may have trouble 
proving that $F \in \cE$!

This completes our sketch of proof of existence, for this specific class $\cE$.
The proof would also work for the nearly Euclidean elliptic integrands of \eqref{3.4}-\eqref{3.5},
but (as far as the author knows) not with more general integrands, because we do not know 
a variant of J. Taylor's theorem in that context.

It seems a little unfortunate that we need to reduce to $d=2$ just for the sake of
finding a local Lipschitz retraction onto $F_k^\ast$. Yet it makes sense to the author
that, especially in the context of minimizing in classes of deformations, the local regularity
of the expected minimizers will play a role in the proof of existence. 
A priori, the existence of retractions looks like a weak regularity property, but the
author is not entirely convinced that it is so weak, and anyway does not know how to
prove it in higher dimensions. This is also why he does not expect an existence result
for the sliding Plateau problem mentioned in Subsection \ref{s2.7} before regularity results
at the boundary.

\ms
We conclude this section with some comments on how the argument above
works with different problems, related to other choices of classes $\cE$.
What follows is more a free discussion than anything, and in particular the author
did not have enough energy to check all the assertions below carefully, but he 
included them anyway because he believes that they explain the sort of issues 
that arise in variants of Plateau problems.

First of all, let us not try to apply the method above when we work with a class $\cE$
which is not stable under deformations (as in Definition \ref{d2.14}). Hopefully this does
not remove too many interesting problems. 

Let us rapidly discuss Reifenberg homology minimizers. The reason why Reifenberg,
and then later authors, restricted to \v{C}ech homology, is that this is the one
that passes to the limit well. Because of that, if we know that the cores $F_k^\ast$
lien in the desired class $\cE$ (thus defined in terms of \v{C}ech homology),
we can eventually get that $F_k^\ast \in \cE$ too. In some cases,
it can be shown that the lower-dimensional part $F_k \sm F_k^\ast$ is useless
for getting the ($d$-dimensional) topological constraints in the definition of $\cE$,
so $F_k^\ast \in \cE$ and we are in business.
If I understand well, this is what happens in \cite{Ita2}, although in the middle
of a very different (and yet beautiful) proof, with the slight disadvantage that the 
authors need to use homology with a compact group (as Reifenberg did), to be
able to cut the $(d-1)$-dimensional piece. The other option, used by Fang \cite{FangEx},
is to keep $F_k \sm F_k^\ast$, make its $(d-1)$-dimensional part converge as an 
Ahlfors-regular set of lower dimension away from $F_k^\ast$, continue with lower
dimensions, and this way get a limit $F$ that is not too big. Life is easier than with
the $d$-dimensional piece, because we don't need to keep precise track of the 
Hausdorff measure of the lower dimensional pieces, just make sure that it stays controlled
so that the $\H^d$ measure of the limit vanishes.

In the presentation above, we carefully avoided the existence of a boundary set $\Gamma$,
because this simplified the discussion. For the Reifenberg homology problem, we cannot do that.
Yet, due to the fact that we only want to minimize the measure $\H^d(E \sm \Gamma)$
away from $\Gamma$, it is possible to limit the discussion about Feuvrier grids and
Federer-Fleming projections to what happens on the complement of $\Gamma$, because the
contribution to a small neighborhood of $\Gamma$ is as small as we want, and for \eqref{9.10}
a control in compact sets of $\cM \sm \Gamma$ is enough. 
See \cite{FangEx} for details.

In contrast, existence results for sliding Plateau problems should probably involve a more careful
study near $\Gamma$, and so would the Reifenberg homology problem if we started
to use a boundary $\Gamma$ with positive $\H^d$-measure, and not require that
$E \supset \Gamma$ in the definitions.

See \cite{LiangTopo} for another existence result based on the Feuvrier scheme, 
\cite{FangHolder} for a case where an existence result (this time for Reifenberg singular
homology minimizers) is deduced from the regularity properties of sliding minimal sets,
and \cite{FangC1} for a similar result for sliding minimal sets of dimension $2$ bounded
by a smooth surface in $\R^3$.

%
%
%
%
%
%
%

\bibspread

\end{document}